
  \magnification 1200


  \newcount\fontset
  \fontset=1
  \def \dualfont#1#2#3{\font#1=\ifnum\fontset=1 #2\else#3\fi}

  \dualfont\bbfive{bbm5}{cmbx5}
  \dualfont\bbseven{bbm7}{cmbx7}
  \dualfont\bbten{bbm10}{cmbx10}

  \font \eightbf = cmbx8
  \font \eighti = cmmi8 \skewchar \eighti = '177
  \font \eightit = cmti8
  \font \eightrm = cmr8
  \font \eightsl = cmsl8
  \font \eightsy = cmsy8 \skewchar \eightsy = '60
  \font \eighttt = cmtt8 \hyphenchar\eighttt = -1
  \font \msbm = msbm10
  
  \font \sixi = cmmi6 \skewchar \sixi = '177
  \font \sixrm = cmr6
  \font \sixsy = cmsy6 \skewchar \sixsy = '60
  \font \tensc = cmcsc10

  \scriptfont \bffam = \bbseven
  \scriptscriptfont \bffam = \bbfive
  \textfont \bffam = \bbten

  \newskip \ttglue

  \def \eightpoint {\def \rm {\fam0 \eightrm }%
  \textfont0 = \eightrm
  \scriptfont0 = \sixrm \scriptscriptfont0 = \fiverm
  \textfont1 = \eighti
  \scriptfont1 = \sixi \scriptscriptfont1 = \fivei
  \textfont2 = \eightsy
  \scriptfont2 = \sixsy \scriptscriptfont2 = \fivesy
  \textfont3 = \tenex
  \scriptfont3 = \tenex \scriptscriptfont3 = \tenex
  \def \it {\fam \itfam \eightit }%
  \textfont \itfam = \eightit
  \def \sl {\fam \slfam \eightsl }%
  \textfont \slfam = \eightsl
  \def \bf {\fam \bffam \eightbf }%
  \textfont \bffam = \bbseven
  \scriptfont \bffam = \bbfive
  \scriptscriptfont \bffam = \bbfive
  \def \tt {\fam \ttfam \eighttt }%
  \textfont \ttfam = \eighttt
  \tt \ttglue = .5em plus.25em minus.15em
  \normalbaselineskip = 9pt
  \def \MF {{\manual opqr}\-{\manual stuq}}%
  \let \sc = \sixrm
  \let \big = \eightbig
  \setbox \strutbox = \hbox {\vrule height7pt depth2pt width0pt}%
  \normalbaselines \rm }



  \newcount \secno \secno = 0
  \newcount \stno \stno =0
  \newcount \eqcntr \eqcntr=0

  \def \ifn #1{\expandafter \ifx \csname #1\endcsname \relax }

  \def \track #1#2#3{\ifn{#1}\else {\tt\ [#2 \string #3] }\fi}

  \def \laberr#1#2{\message{*** RELABEL CHECKED FALSE for #1 ***}
      RELABEL CHECKED FALSE FOR #1, EXITING.
      \end}

  \def \seqnumbering {\global \advance \stno by 1 \global
    \eqcntr=0 \number \secno .\number \stno }

  \def \current {\number \secno
    \ifnum \number \stno = 0\else .\number \stno \fi }

  \def \eqmark#1{\global \advance\eqcntr by 1
    \edef\a{\number \secno .\number\stno.\number\eqcntr}
    \eqno {(\a)}
    \global \edef #1{\a}\track{showlabel}{*}{#1}}

  \def \syslabel#1#2{\global \expandafter \edef \csname
    #1\endcsname {#2}}

  \def \fcite#1#2{\syslabel{#1}{#2}\lcite{#2}}

  \def \label #1 {%
    \ifn {#1}%
      \syslabel{#1}{\current}%
    \else
      \edef\a{\expandafter\csname #1\endcsname}%
      \edef\b{\current}%
      \ifx \a \b \else \laberr{#1=(\a)=(\b)} \fi
      \fi
    \track{showlabel}{*}{#1}}

  \def \lcite #1{(#1\track{showcit}{$\bullet$}{#1})}

  \def \cite #1{[{\bf #1}\track{showref}{\#}{#1}]}

  \def \scite #1#2{{\rm [\bf #1\track{showref}{\#}{#1}{\rm \hskip 0.7pt:\hskip 2pt #2}\rm]}}


 \def \Headlines #1#2{\nopagenumbers
    \advance \voffset by 2\baselineskip
    \advance \vsize by -\voffset
    \headline {\ifnum \pageno = 1 \hfil
    \else \ifodd \pageno \tensc \hfil \lcase {#1} \hfil \folio
    \else \tensc \folio \hfil \lcase {#2} \hfil
    \fi \fi }}

  \def \Date #1 {\footnote {}{\eightit Date: #1.}}


  \def \lcase #1{\edef \auxvar {\lowercase {#1}}\auxvar }

  \def \goodbreak {\vskip0pt plus.1\vsize \penalty -250 \vskip0pt
plus-.1\vsize }

  \def \section #1{\global\def \SectionName{#1}\stno = 0 \global
\advance \secno by 1 \bigskip \bigskip \goodbreak \noindent {\bf
\number \secno .\enspace #1.}\medskip \noindent \ignorespaces}

  \long \def \sysstate #1#2#3{\medbreak \noindent {\bf \seqnumbering
.\enspace #1.\enspace }{#2#3\vskip 0pt}\medbreak }
  \def \state #1 #2\par {\sysstate {#1}{\sl }{#2}}
  \def \definition #1\par {\sysstate {Definition}{\rm }{#1}}
  \def \remark #1\par {\sysstate {Remark}{\rm }{#1}}


  \def \proof {\medbreak \noindent {\it Proof.\enspace }}
  \def \proofend {\ifmmode \eqno \square \else \hfill \square
\looseness = -1 \medbreak \fi }

  \def \$#1{#1 $$$$ #1}
  \def \=#1{\buildrel #1 \over =}

  \def \Item #1{\smallskip \item {#1}}
  \newcount \zitemno \zitemno = 0
  \def \izitem {\zitemno = 0}
  \def \zitem {\global \advance \zitemno by 1 \Item {{\rm(\romannumeral
\zitemno)}}}

  \newcount \footno \footno = 1
  \newcount \halffootno \footno = 1
  \def \footcntr {\global \advance \footno by 1
  \halffootno =\footno
  \divide \halffootno by 2
  $^{\number\halffootno}$}
  \def \fn#1{\footnote{\footcntr}{\eightpoint#1}}


  \def \N {{\bf N}}
  \def \C {{\bf C}}
  \def \<{\left \langle \vrule width 0pt depth 0pt height 8pt }
  \def \>{\right \rangle }  
  \def \curly#1{{\cal #1}}
  
  \def \and {\hbox {,\quad and \quad }}
  \def \calcat #1{\,{\vrule height8pt depth4pt}_{\,#1}}
  \def \labelarrow#1{\ {\buildrel #1 \over \longrightarrow}\ }
  \def \for #1{,\quad \forall\,#1}
  \def \square {\hbox {$\sqcap \!\!\!\!\sqcup $}}
  \def \crossproduct {{\hbox {\msbm o}}}
  \def \stress #1{{\it #1}\/}
  \def \inv {^{-1}}
  \def \*{\otimes}

  \newcount \bibno \bibno =0
  \def \newbib #1{\global\advance\bibno by 1 \edef #1{\number\bibno}}
  \def \bibitem #1#2#3#4{\smallskip \item {[#1]} #2, ``#3'', #4.}
  \def \references {
    \begingroup
    \bigskip \bigskip \goodbreak
    \eightpoint
    \centerline {\tensc References}
    \nobreak \medskip \frenchspacing }


  \def \O{{\cal O}}

  \def \Integers{{\bf Z}}
  \def \compos{\circ}
  \def \({\left(}
  \def \){\right)}
  \def \tilde#1{\widetilde {#1}}
  \def \phi{\varphi}
  \def \ext#1#2{\varepsilon_{_{#1#2}}}
  \def \implies{\ \Rightarrow\ }

  \def \Kone{{\cal K}^0}

  \def \K{{\cal K}}
  \def \M{{\cal M}}
  \def \H{{\cal H}}
  \def \Z{{\cal Z}}

  \def \tK{\widehat{\cal K}}
  \def \tM{\widehat{\cal M}}
  \def \tZ{\widehat{\cal Z}}

  \def \LL#1#2{\curly {L}_{#1}(\H_{#2})}
  \def \spr{\mu}

  \font\forsmallbullet=cmsy10 scaled 300
  \def \smallbullet{\hbox{\forsmallbullet }}

  \def \prd#1{{\buildrel \smallbullet \over #1}}
  \def \wprd#1{(#1)^{\raise 3pt \hbox{\forsmallbullet }}}

  \def \ts{\widehat s}
  \def \tidp#1{\widehat e_{#1}}
  \def \tj{\widehat{\hbox{\j}}}  
  \def \rxv#1#2{#1{\times}#2}
  \def \trxv#1#2{\widehat{\rxv{#1}{#2}}}
  \def \CpPi{\rxv i\n}
  \def \TPi{\widehat\CpPi}   \def \TPi{\trxv i\n}
  \def \cpj{\hbox{\j}}
  \def \cps{s}
  \def \m{V}
  \def \ac{\alpha}
  \def \n{\check s}
  \def \e{\varepsilon}  \def \e{E}
  \def \R{{\cal R}}
  \def \L{{\cal L}}
  \def \l{\ell}
  \def \I{{\cal I}}
  \def \D{{\cal D}}
  \def \B{{\cal B}}
  \def \CB{C^*(\B)}
  \def \*{\otimes}
  \def \Lin{\curly {L}}
  \def \toep{\curly{T}(A,G,\m)}
  \def \syscp#1{A \crossproduct_{#1} G}
  \def \sgcp{A\crossproduct_\ac \Pos}
  \def \cp{\syscp{{_\m}}}
  
  \def \pe#1{\vrule depth #1pt width 0pt}
  \def \inc#1{#1}
  \def \Pos{P}

  \def \letrep{v}
  \def \rep#1{\letrep_{#1}}
  \def \idp#1{e_{#1}}
  \def \wd{\alpha} 
  \def \wdaux{\beta}

  \def \bratteli{BRo} \newbib\bratteli
  \def \brownlowe{BRa} \newbib\brownlowe
  \def \cuntz{C} \newbib\cuntz
  \def \cuntzTwo{C2} \newbib\cuntzTwo
  \def \amena{E1} \newbib\amena
  \def \inverse{E2} \newbib\inverse
  \def \endo{E3} \newbib\endo
  \def \interaction{E4} \newbib\interaction
  \def \vershik{EV} \newbib\vershik
  \def \fell{FD} \newbib\fell
  \def \fowler{F} \newbib\fowler
  \def \katsura{K} \newbib\katsura
  \def \wasserman{W} \newbib\wasserman
  \def \lacareaburn{LR} \newbib\lacareaburn
  \def \larsen{L} \newbib\larsen
  \def \cpim{P} \newbib\cpim
  \def \renault{R} \newbib\renault
  \def \takesaki{T} \newbib\takesaki

  \centerline{\bf A NEW LOOK AT THE CROSSED-PRODUCT OF A}
  \smallskip
  \centerline{\bf  C*-ALGEBRA BY A SEMIGROUP OF}
  \smallskip
  \centerline{\bf ENDOMORPHISMS}

  \Headlines
  {Crossed products by endomorphisms}
  {R.~Exel}

  \null

  \bigskip
  \centerline{{\tensc Ruy Exel}\footnote{$^\star$}{\eightpoint
Partially supported by CNPq.}}
  \bigskip

  \Date{17 Nov 2005}

  \midinsert 
  \narrower \narrower
  \eightpoint \noindent {\tensc Abstract}.
  Let $G$ be a group and let $\Pos\subseteq G$ be a subsemigroup.  In
order to describe the crossed product of a C*-algebra $A$ by an action
of $\Pos$ by unital endomorphisms we find that we must extend the
action to the whole group $G$.  This extension fits into a broader
notion of \stress{interaction groups} which consists of an assignment
of a positive operator $\m_g$ on $A$ for each $g$ in $G$, obeying a
partial group law, and such that $(\m_g,\m_{g\inv})$ is an interaction
for every $g$, as defined in a previous paper by the author.  We then
develop a theory of crossed products by interaction groups and compare
it to other endomorphism crossed product constructions.
  \endinsert

  \section{Introduction}
  Roughly five years ago I wrote a similarly titled paper \cite{\endo}
in which a new notion of crossed product by an endomorphism was
introduced.  One of the reasons for doing so was a certain
dissatisfaction with the existing crossed product theories because
none of these could give the ``correct'' answer in the example of
Markov subshifts.

I explain. Let $A$ be an $n\times n$ matrix with 0--1 entries and
$K\subseteq\{1,\ldots,n\}^{\N}$ be the corresponding Markov space.
Denoting by $T$ the shift on $K$ given by
  $T(x)_n = x_{n+1}$,
  one may define an endomorphism $\ac$ on the C*-algebra $C(K)$ by
setting $\ac(f) = f\compos T$, for all $f$ in $C(K)$. 

It has always been my impression that, whatever crossed product
construction one could possibly conceive of, it should give the
Cuntz--Krieger algebra $\O_A$ if applied to the endomorphism
described above.  This is of course a matter of personal opinion, but
one must concede that the close relationship between Markov subshifts
and Cuntz-Krieger algebras definitely suggests this.
  The construction described in \cite{\endo} is therefore the result of
my search for a theory of crossed products which fulfills these
expectations.


The basic starting point in the construction of any C*-algebra
associated to a given endomorphism $$ \ac:A\to A$$ is to decide what
is one's idea of \stress{covariant representation}.  In other words
one imagines that $A$ is a subalgebra of some bigger C*-algebra where
the given endomorphism can be implemented by means of some algebraic
operation.

A very popular proposal is to require the \stress{covariance condition}
  $$
  \ac(a) = SaS^*,
  \eqno{(\seqnumbering)}
  \label FirstCovariance
  $$
  where $S$ is and isometry in this bigger universe, that is, an
element satisfying $S^*S=1$.  This expression is very appealing
because, after all, the correspondence
  $
  a\mapsto SaS^*
  $
  is a *-homomorphism for any isometry $S$!

However, even though the algebra of continuous functions on Markov's
space is a subalgebra of $\O_A$ in a standard way, i.e.~as the
diagonal of the standard AF-subalgebra, there is no isometry $S$ in
$\O_A$ which implements the Markov subshift according to the above
covariance condition.  Please see the Appendix for a detailed proof of
this statement.

 When one studies nonunital endomorphisms,
especially those with hereditary range, condition
\lcite{\FirstCovariance} becomes a powerful tool leading to very
interesting applications, but unfortunately it must be abandoned if
one deals with unital endomorphisms.

Carefully studying the case of the Cuntz-Krieger algebra I eventually
settled for the following substitute of \lcite{\FirstCovariance}:
  $$
  \matrix{ Sa=\ac(a)S, \pe{6}\cr
  S^*aS = \L(a).}
  $$
  While the first condition has already been discussed in the
literature, the appearance of the \stress{transfer operator} $\L$ was
a new feature.  It should be stressed that $\L$ needs to be given in
advance; it is an important ingredient of the construction inasmuch as
the endomorphism $\ac$ itself.  See \cite{\endo} for more details.

The resulting theory of crossed products based on transfer operators
seems to be getting more and more acceptance and hence it is natural
to try to apply the same ideas for semigroups of endomorphisms.
Suppose, therefore, that $A$ is a unital C*-algebra, $\Pos$ is a
semigroup and 
  $$
  \ac:\Pos \to {\rm End}(A),
  $$
  is a semigroup homomorphism, a.k.a.~an action by endomorphisms.

In case the range of every $\ac_g$ is a hereditary subalgebra of $A$
the natural covariance condition is definitely to ask for a semigroup
of isometries $\{S_g\}_{g\in \Pos}$ such that
  $$
  \ac_g(a) = S_g aS_g^*,
  $$
  as done e.g.~in \scite{\lacareaburn}{2.2}. But, as seen above, this
is not appropriate for unital endomorphisms.

Given that we are now dealing with many endomorphisms at the same
time, we should also need many transfer operators.  However one must
decide how to put these together in a harmonious way.

If one resorts to the single endomorphism case for guidance, given
$\ac$ and $\L$, it is easy to see that $\L^n$ is a transfer operator
for $\ac^n$.  It is also interesting to notice that one then gets a
collection $\{\m_n\}_{n\in\Integers}$ of bounded linear maps on $A$
given by
  $$
  \m_n =
  \left\{\matrix{ \ac^n, & \hbox{if } n\geq0, \cr\cr
                   \L^{-n}, & \hbox{if } n<0.}\right.
  \eqno{(\seqnumbering)}   
  \label SimpleInteracGroup
  $$
  Thus, not only the positive integers are assigned to an operator on
$A$ (namely $\ac^n$), but the negative integers are similarly given a
role.  This is related to the discussion in the introduction of
\cite{\interaction} regarding past time evolution for irreversible
dynamical systems.

  Assuming that $\Pos$ is a subsemigroup of a group $G$ it is then
tempting to generalize by requiring that one is given an operator
$\L_h$ for each $h\in\Pos\inv$, such that $\L_{g\inv}$ is a transfer
operator for $\ac_g$, for every $g$ in $\Pos$.  This is essentially
what is done in \cite{\larsen}.

However it seems a bit awkward to me that elements of
$\Pos\cup\Pos\inv$ are allowed to act on $A$, i.e., are assigned an
operator on $A$ (either an endomorphism or a transfer operator) but
the other elements in the group are not given any role.  On totally
ordered groups, that is, when $G=\Pos\cup\Pos\inv$, there are really
no other such elements but one may be interested in more general
group-subsemigroup pairs such as $(\Integers^2,\N^2)$.

We therefore face a crucial question.  Given $g$ not in
$\Pos\cup\Pos\inv$, should we assign to it an endomorphism or a
transfer operator?  It was precisely my interest in solving this
puzzle that led me to introduce the notion of \stress{interactions}
\cite{\interaction}.  These are operators on C*-algebras which
simultaneously generalize endomorphisms and transfer operators.

To cut a long story short I was led to reformulating the very notion
of semigroups of endomorphisms: it should be replaced by an
assignment of an operator $\m_g$ for every $g$ in $G$, such that the
pair $(\m_g,\m_{g\inv})$ is an interaction.

It could perhaps be the case that $G$ contains a subsemigroup $\Pos$
such that $\m_g$ is an endomorphism for every $g$ in $\Pos$, in which
case $\m_{g\inv}$ will necessarily be a transfer operator for $\m_g$.
But then again this may not be so.  From this point on I decided not
to rely on the existence of subsemigroups any more and fortunately the
price to be payed was not so high.

Of course we have so far neglected to consider which kind of
dependence the maps $\m_g$ should have on the variable $g$.  Obviously
it would be too much to require that $\m_g\m_h$ coincide with
$\m_{gh}$, as this is not even true for \lcite{\SimpleInteracGroup}.
On the other hand it would be crazy not to require any compatibility
between the group law and the composition of the $\m_g$.

Resorting again to the single endomorphism case I was  able to prove
that the $\m$ of \lcite{\SimpleInteracGroup} is a \stress{partial
representation} of $\Integers$ on $A$, meaning that
  $$
  \matrix{ \m_n\m_m\m_{-m} = \m_{n+m}\m_{-m} \pe{9}\cr
  \m_{-n}\m_n\m_m = \m_{-n}\m_{n+m}.}
  $$
  This says that the group law
  ``$\m_n\m_m = \m_{n+m}$'' does hold, as long as it is
right-multiplied by $\m_{-m}$, or left-multiplied by $\m_{-n}$.

Following this I was led to the definition of what we shall call an
\stress{interaction group}: it is a triple $(A,G,\m)$, where $A$ is a
C*-algebra, $G$ is a group, and $\m$ is a partial representation of
$G$ on $A$ such that $(\m_g,\m_{g\inv})$ is an interaction for every
$g$ in $G$.

I am well aware that the search for a \stress{new definition} is
always a risky business unless one has valuable examples under one's
belt.  In fact it was precisely this fear that led me to keep the
findings of this paper to myself for quite some time now.  However the
problem of finding an appropriate definition for the notion of crossed
product by semigroups of unital endomorphisms has been around for
almost thirty years now \scite{\cuntz}{\S 2}.  In addition, recent
discussions with Jean Renault suggested the existence of some quite
interesting examples which we hope to describe in a forthcoming paper.

A few words must now be said about the crossed product.  As already
mentioned the starting point should be the covariance condition.
Using the ideas of \cite{\interaction} it is not so hard to arrive at
the following: a covariant representation of an interaction group
$(A,G,\m)$ in a C*-algebra $B$ consists of a *-homomorphism
  $$
  \pi:A\to B
  $$
  and a *-partial representation 
  $$
  \letrep: G\to B,
  $$
  such that 
  $$
  \rep g \pi(a) \rep{g\inv} = \pi\(\m_g(a)\)   \rep g \rep{g\inv},
  $$  
  for every $g$ in $G$, and every $a$ in $A$.
  Once this is accepted the experienced reader will immediately think
of the universal C*-algebra for covariant representations, which we
shall denote $\toep$.

However this is definitely too big an object since among its
representations one finds bizarre things such as $\letrep_g\equiv 0$.
This algebra must therefore be recognized as just a first step in the
construction of the crossed product.  A fine tunning is necessary to
restrict its wild collection of representations.

This fine tuning consists in eliminating the ``redundancies'', as done
in \cite{\endo} and cleverly recognized by Brownlowe  and Raeburn
\cite{\brownlowe} as the passage from the Toeplitz-Cuntz-Pimsner
algebra to the Cuntz-Pimsner algebra in the case of a correspondence
\scite{\cpim}{3.4 and 3.12}.

There are in fact several other instances in the literature where
ideas with a similar flavor appeared.  Among these we cite 
  Larsen's redundancies \scite{\larsen}{1.2} and Fowler's
Cuntz-Pimsner covariance in the context of discrete product systems
\scite{\fowler}{2.5}.

Having barely no examples in hand to help decide what the correct
notion of redundancies should be, I had to search for a clue
somewhere.  Recall that, in the context of a single
endomorphism/transfer operator pair $(\ac,\L)$, one has a very
concrete model of the crossed product, provided one is given a
faithfull invariant state.  In this case $A\crossproduct_{\ac,\L}\N$
essentially\fn{One in fact needs to tensor this with the regular
representation of $\Integers$.} coincides with the algebra of
operators on the GNS space of the given state, generated by $A$ and
the isometry naturally associated to $\ac$ \scite{\vershik}{6.1}.

One nice feature of interaction groups is that, given an invariant
state, one can also easily construct a covariant representation
$(\pi,\letrep)$ on the GNS space of the given state
\fcite{GNSStrongCov}{11.4}.  I therefore thought it would be nice if
the crossed product ended up isomorphic to the concrete algebra of
operators generated by $\pi(A)$ and $\letrep(G)$ when the state is
faithfull, precisely as in the case described above.  Of course the
latter algebra could not be taken as the definition of the crossed
product, unless of course one postulated the existence of a faithfull
invariant state.

I therefore set out to find the right notion of redundancies guided
by the requirement that, once these are moded out, one should arrive at a
quotient of $\toep$ which in turn is isomorphic to the above concrete
algebra of operators provided a faithfull invariant state exists.

This was a rather interesting clue to follow because, on the one hand,
it suddenly put our very abstract theory in close contact with a very
concrete algebra of operators and, on the other, it lead to the highly
complex notion of redundancies given in Definition
  \fcite{DefineRedundancies}{6.1} which I would not have imagined
otherwise.
  I would even go as
far as to suggest that this notion of redundancies could have some
impact in a possible strengthening of the notion of Cuntz-Pimsner
covariance in the context of discrete product systems.

We feel that the questions we address here are complex enough to
warrant avoiding extra difficulties.  For that reason we chose to
simplify things as much as possible by assuming two extra standing
hypothesis, namely that interactions preserve the unit
\fcite{DefGinter}{3.1} and that the conditional expectations
$\m_g\m_{g\inv}$ one gets from an interaction are non-degenerated
\fcite{DefineNonDegInter}{3.3}.  In the same way that the
Cuntz--Pimsner construction \cite{\cpim} was later refined by Katsura
\cite{\katsura} we believe that our results may be similarly extended
to avoid these special assumptions.

The organization of the paper is as follows:  
  in section (2) we give a few basic facts about partial
representations and in section (3) we inroduce the central notion of
interaction groups.
  In section (4) we study covariant representations and collect some
crucial technical tools in preparation for the introduction of the
notion of redundancies.

  The Toeplitz algebra of an interaction group is defined next and in
section (6) the all important notion of redundancies is finally given,
together with the definition of the crossed product.

  In section (7) we show that the crossed product possesses a natural
grading over the group $G$.  

Up to this point the development is admitedly abstract but in section
(8) we study at lenght certain Hilbert modules and operators between
them to be used in the following section with the purpose of finally
giving a concrete representation of the crossed product.  This may be
thought of as a generalization of the notion of regular
representation.  As a byproduct we show that the natural map from $A$
to $\cp$ is injective.

In section (10) we study the general problem of faithfulness of
representations of the crossed product culminating with Theorem
\fcite{RegularIsFaithful}{10.10} in which we show that the regular
representation of the crossed product is faithful provided $G$ is
amenable.

In section (11) we assume the existence of a faithful invariant state
and show in Theorem \fcite{ConcreteCp}{11.7} the result already
alluded to above, according to which $\cp$ is isomorphic to a concrete
algebra of operators on a certain amplification of the GNS space,
again under the assumption that $G$ is amenable.
  We hope that this Theorem will convince the skeptical reader that
our theory is not as outlandish as it might seem.

In section (12) we finally return to considering semigroups of
endomorphisms as well as the question of whether such a semigroup can
be extended to give an interaction group.  Assuming the existence of
certain states we show in Theorem \fcite{UniqueExtAndCp}{12.3} that,
when such an extension exists, it is unique and that $\cp$ is
isomorphic to a natural concretely defined endomorphism crossed
product.  This could support the argument that interaction groups
underlie certain semigroups of endomorphisms.

In the following section we compare our theory with \cite{\larsen}
giving a necessary and sufficient condition for Larsen's dynamical
systems to be extended to an interaction group.  In section \lcite{14}
we show, by means of an example, that Larsen's systems do not always
extend to an interaction group.

  Given a substantial difference in the respective notions of
redundancies we do not believe that our crossed product coincides with
Larsen's.  However it might be that these will coincide after Larsen's
redundancies are replaced by some more general notion involving
several semigroup elements at the same time.

I would like to acknowledge some very interesting conversations with
Jean Renault, based on which I gathered enough courage to finish a
rough draft of this work which was laying dormant for roughly three
years now.  I would also like to express my thanks to Vaughan Jones
for having suggested the very interesting example described in section
\lcite{14}.  Thanks go also to 
  Joachim Cuntz,
  Jean Renault,
  Mikael R\o rdam, and
  Klaus Thomsen,
  for helping me to sort out relevant bibliography.

\section{Partial representations}
  In this section we discuss some elementary facts about partial
representations which will be used in the sequell.
  The reader is referred to \cite{\inverse} for more information.

\definition 
  \label DefPrep
  A \stress{partial representation} of a group $G$ in a
unital algebra $A$ is a map
  $$
  \letrep : G \to A
  $$
  such that for all $g$ and $h$ in $G$ one has that
  \izitem
  \zitem $\rep 1 = 1$,
  \zitem $\rep g\rep h\rep {h\inv} = \rep {gh}\rep {h\inv}$, and
  \zitem $\rep {g\inv}\rep g\rep h = \rep {g\inv}\rep {gh}$.

The following is a usefull basic result about partial representations:

\state Proposition 
  \label BasicPreps
  Given a partial representation $\letrep$ let
  $$
  \idp g = \rep g\rep {g\inv}.
  $$
  Then for all $g$ and $h$ in $G$,
  \izitem 
  \zitem $\idp g$ is an idempotent,
  \zitem $\rep g \idp h = \idp {gh} \rep g$,
  \zitem $\idp g$ commutes with $\idp h$.

  \proof See \scite{\inverse}{2.4}.
  \proofend

Let us now study one-sided invertible elements in partial representations:

\state Lemma
  \label OnOneSidedInverse
  Fix a partial representation $\letrep$.
  \izitem 
  \zitem Given $g$ in $G$ suppose that $u$ is a right-inverse of $\rep
g$, in the sense that $\rep g u =1$.  Then $\rep{g\inv}$ is also a
right-inverse of $\rep g$.
  \zitem If $\rep g$ is left-invertible then $\rep{g\inv}$ is  a
left-inverse of $\rep g$.
  \zitem If $\rep g$ is right-invertible then $\rep g \rep h =
\rep{gh}$, for all $h$ in $G$.
  \zitem If $\rep g$ is left-invertible then $\rep h \rep g =
\rep{hg}$, for all $h$ in $G$.

\proof
  We have 
  $$
  \rep g \rep {g\inv} = 
  \rep g \rep {g\inv} \rep g u= 
  \rep g u=  1.
  $$
  This proves (i) while (ii) may be proved in a similar way.  As for
(iii) we have by (i) that
  $$
  \rep g \rep h = 
  \rep g\rep{g\inv} \rep g \rep h =
  \rep g\rep{g\inv} \rep {gh} =
  \rep {gh}.
  $$
  The last point follows similarly.
  \proofend

The following is often usefull when one deals with ordered groups:

  \state Lemma 
  \label RepsCoincidingOnPos
  Let $\Pos$ be a subsemigroup of $G$ such that $G=\Pos\inv\Pos$, and
let $\letrep$ be a partial representation of $G$ in $A$ such that
$\rep g$ is left-invertible for every $g$ in $\Pos$, then
  $$
  \rep{x\inv y}= \rep {x\inv} \rep y
  \for x,y\in\Pos.
  $$
  Moreover, if $\letrep'$ is another partial representation of $G$ in
$A$ such that $\rep g' = \rep g$ for all $g\in \Pos\cup\Pos\inv$, then
$\letrep'=\letrep$.

  \proof That $\rep{x\inv y}= \rep {x\inv} \rep y$ 
follows from \lcite{\OnOneSidedInverse.iv}.  That $\letrep'=\letrep$
is then obvious.
  \proofend

By a \stress{word} in $G$ we will mean a finite sequence
  $$
  \wd = (g_1,g_2,\ldots,g_n)
  \eqno{(\seqnumbering)}
  \label DefineWord
  $$
  of elements in $G$.  Given a word $\wd$ as above we will let
  $$
  \wd\inv = (g_n\inv,g_{n-1}\inv,\ldots,g_1\inv).
  $$

Fixing, for the time being, a partial representation $\letrep$ of
$G$ on the algebra $A$, we will let for every word $\wd$ as above,
  $$
  \rep \wd = \rep {g_1}\rep {g_2}\ldots\rep {g_n},
  \eqno{(\seqnumbering)}
  \label DefineRepOnWords
  $$ and $$
  \idp \wd = \rep \wd \rep {\wd\inv}.
  $$

\state Proposition
  \label BasicPrepForWords
  Let $\wd = (g_1,g_2,\ldots,g_n)$ be a word in $G$.   Then
  \izitem
  \zitem $\idp \wd \rep \wd = \rep \wd$.
  \zitem $\idp \wd = \idp {g_1}\idp {g_1g_2}\idp {g_1g_2g_3} \ldots
\idp {g_1g_2g_3\ldots g_n}$, and hence $\idp \wd$ is  idempotent,
  \zitem $\idp \wd = \idp {g_1g_2\ldots g_n} \idp \wdaux$, where
$\wdaux=(g_1,g_2,\ldots,g_{n-1})$.

  \proof
  In order to prove (iii) observe that
  $$
  \idp \wd =
  \rep \wdaux \rep{g_n}\rep{g_n\inv} \rep{\wdaux\inv} =
  \rep \wdaux \idp {g_n} \rep{\wdaux\inv} = \cdots
  $$
  Applying \lcite{\BasicPreps.ii} repeatedly we conclude that the
above equals
  $$
  \cdots = 
  \idp {g_1g_2\ldots g_n}   \rep \wdaux  \rep{\wdaux\inv} =
  \idp {g_1g_2\ldots g_n}   \idp \wdaux,
  $$
  thus proving (iii).  It is now easy to show that (ii) follows from
(iii) and by induction.
  As for (i) we have by (iii) that
  $$
  \idp \wd \rep \wd = 
  \idp {g_1g_2\ldots g_n}   \idp \wdaux  \rep\wdaux \rep {g_n} =\cdots
  $$
  By induction the above equals
  $$
  \cdots =
  \idp {g_1g_2\ldots g_n} \rep\wdaux \rep {g_n} =
  \rep\wdaux \idp {g_n} \rep {g_n} =
  \rep\wdaux \rep {g_n} \rep {g_n\inv}\rep {g_n} \={(\DefPrep.ii)}
  \rep\wdaux \rep {g_n} =
  \rep\wd.
  \proofend
  $$

  If $\wd = (g_1,g_2,\ldots,g_n)$ is a word in $G$ we will denote by
$\prd \wd$ the  product of all components of $G$, namely
  $$
  \prd \wd = g_1g_2g_3\ldots g_n.
  \eqno{(\seqnumbering)}
  \label DefinePrd
  $$
We will also let 
$\spr(\wd)$ be the subset of $G$ given by
  $$
  \spr(\wd) = \{1,g_1,g_1g_2,g_1g_2g_3, \ldots, g_1g_2g_3\ldots g_n\},
  \eqno{(\seqnumbering)}
  \label DefMuWd
  $$
  so that $\spr(\wd)$ consists of the $\prd \wdaux$ for all initial
segments $\wdaux$ of $\wd$.  Notice that we have included $\prd
\wdaux$ for the empty initial segment $\wdaux$, namely 1.

  Employing this notation we observe that
\lcite{\BasicPrepForWords.ii} can be stated as
  $
  \idp \wd = \prod_{h\in\spr(\wd)}\idp h,
  $
  while \lcite{\BasicPrepForWords.iii} reads
  $
  \idp \wd = \idp {\prd \wd} \idp \wdaux.
  $

The following elementary property of words will be usefull later on:

  \state Proposition
  \label SProdConcat
  Given words $\wd$ and $\wdaux$ in $G$ denote by $\wd\wdaux$ the
concateneted word.  Then
  $$
  \spr(\wd\wdaux) = \spr(\wd) \cup \prd \wd \spr (\wdaux).
  $$

Another curious fact which will become relevant  is:

\state Proposition
  \label SprOfInverse
  Let $\wd$ be a word in $G$.  Then
  \izitem
  \zitem $\spr(\wd) = \prd\wd\spr(\wd\inv)$,
  \zitem If $\prd\wd=1$, then $\spr(\wd)=\spr(\wd\inv)$.

  \proof
  Let $\wd=(g_1,\ldots,g_n)$.
  Since  
  $$
  g_1g_2\ldots g_{k-1}g_k \ldots g_{n-1}g_n=\prd\wd,
  $$
  we have that
  $$
  g_1g_2\ldots g_{k-1} = \prd\wd g_n\inv g_{n-1}\inv\ldots g_k\inv
  \for k=1,\ldots,n,
  $$
  from which (i) follows.  Obviously (i) implies (ii).
  \proofend

  If $\wd=(g_1,\ldots,g_n)$ is a word in $G$ recall that
  $$
  \rep \wd = \rep {g_1}\rep {g_2}\ldots\rep {g_n},
  $$
  while 
  $$
  \rep{\prd\wd} = \rep {g_1g_2\ldots g_n}.
  $$
  In our next Proposition we will relate these.

  \state Proposition
  \label RepWdVsRepPrdWd
  For every word $\wd$ in $G$ one has $\rep\wd=\idp\wd\rep{\prd\wd}$

  \proof
  Let  $\wd=(g_1,\ldots,g_n)$.  If $n=1$ the statement follows
immediately from \lcite{\DefPrep.ii} with $g=g_1$ and $h=g_1\inv$.
Proceeding by induction on $n$ let $\wdaux=(g_2,\ldots,g_n)$ and observe that
  $$
  \rep\wd =
  \rep{g_1} \rep\wdaux =
  \rep{g_1} \idp\wdaux \rep{\prd\wdaux} =
  \rep{g_1}\rep{g_1\inv}  \rep{g_1} \idp\wdaux \rep{\prd\wdaux} \={(\BasicPreps.iii)}
  \rep{g_1}\idp\wdaux \rep{g_1\inv}  \rep{g_1} \rep{\prd\wdaux}
\={(\DefPrep.iii)}
  $$$$=
  \rep{g_1}\rep\wdaux\rep{\wdaux\inv} \rep{g_1\inv}  \rep{\prd\wd} =
  \rep\wd\rep{\wd\inv}  \rep{\prd\wd} =
  \idp\wd \rep{\prd\wd}
  \proofend
  $$

  \state Proposition
  \label ComparisonEalfaEbeta
  If $\wd$ and $\wdaux$ are words in $G$ such that
$\spr(\wd)\subseteq\spr(\wdaux)$ then $\idp\wd\geq\idp\wdaux$ in the
usual sense that
  $\idp\wd\idp\wdaux=\idp\wdaux$.

  \proof
  Follows immediately from \lcite{\BasicPrepForWords.ii}.
  \proofend

  The following simple result will be of crucial importance later on:

  \state Corollary
  \label TheFactIMissed
  Let $\wd$ and $\wdaux$ be words in $G$ such that $\prd\wd=1$ and
$\spr(\wd)\subseteq\spr(\wdaux)$.  Then $\rep\wd\rep\wdaux=\rep\wdaux$.

  \proof
  Since $\rep{\prd\wd}=\rep1=1$ we have by \lcite{\RepWdVsRepPrdWd}
that $\rep\wd=\idp\wd$.  Thus
  $$
  \rep\wd\rep\wdaux  =
  \idp\wd\rep\wdaux  \={(\BasicPrepForWords.i)}
  \idp\wd\idp\wdaux\rep\wdaux  \={(\ComparisonEalfaEbeta)}
  \idp\wdaux\rep\wdaux  =
  \rep\wdaux.
  \proofend
  $$

  \definition If $A$ is a C*-algebra and $\letrep:G\to A$ is a partial
representation we will say that $\letrep$ is a *-partial
representation if 
  $$
  \rep g ^* =
  \rep {g\inv}
  \for g\in G.
  $$

\state Proposition
  If $\letrep$ is a *-partial representation then for every word $\wd$
in $G$ one has that
  $\rep \wd$ is a partial isometry and 
  $\idp \wd$ is a projection (self-adjoint idempotent).

  \proof It is clear that $\rep\wd^* = \rep{\wd\inv}$.  So 
  $\idp \wd = \rep \wd \rep \wd^*$ is clearly self-adjoint. 
  Moreover
  $$
  \rep \wd  \rep \wd^*  \rep \wd =
  \idp \wd  \rep \wd \={(\BasicPrepForWords.i)}
  \rep \wd,
  $$
  so $\rep\wd$ is indeed a partial isometry.
  \proofend 

  \section{Interaction groups}
  In this section we introduce the main object of this work.  Given a
Banach space $X$ we denote by ${\cal B}(X)$ be the algebra of all
bounded linear operators on $X$.  Below we will refer to ${\cal
B}(A)$, where $A$ is a C*-algebra.

  \definition
  \label DefGinter
  An \stress{interaction group} is a triple $(A,G,\m)$ such that
$A$ is a unital C*-algebra, $G$ is a group, and 
  $$
  \m : G \to {\cal B}(A)
  $$
  is a partial representation such that, for every $g$ in $G$,
  \izitem
  \zitem $\m_g$ is a positive map,
  \zitem $\m_g(1)=1$, 
  \zitem for every $a$ and $b$ in $A$, such that either $a$ or $b$
belongs to the range of $\m_{g\inv}$, one has that
$\m_g(ab)=\m_g(a)\m_g(b)$.

\medskip \noindent
Observe that for each $g$ in $G$, the pair  $(\m_g,\m_{g\inv})$ is an
\stress{interaction} according to \scite{\interaction}{Definition 3.1}.

We will always denote by $\R_g$ the range of $\m_g$, and we will  let
  $$
  \e_g = \m_g\m_{g\inv}.
  $$

The following is a direct consequence of \scite{\interaction}{2.6, 2.7
and 3.3}:

  \state Proposition
  \label CompletelyPos
  Given an interaction-group $(A,G,\m)$, for every $g$ in $G$ one has
that
  \izitem
  \zitem $\m_g$ is completely positive and completely contractive,
  \zitem $\R_g$ is a closed *-subalgebra of $A$,
  \zitem $\e_g$ is a conditional expectation onto $\R_g$,
  \zitem $\m_g$ restricts to a *-isomorphism from $\R_{g\inv}$ onto
$\R_g$, whose inverse is the corresponding restriction of $\m_{g\inv}$.

  Recall that a conditional expectation $E$ is said to be non-degenerated if 
  $$
  E(a^*a) = 0 \implies a=0.
  $$

\definition
  \label DefineNonDegInter
  An interaction-group $(A,G,\m)$ will be said to be non-degenerated
if for every $g\in G$ one has that $\e_g$ is a non-degenerated
conditional expectation.

  From now on, and throughout the rest of this paper, we will fix a
non-degenerated interaction-group $(A,G,\m)$.

If $\wd=(g_1,\ldots,g_n)$ is a word in $G$, 
  in accordance with \lcite{\DefineRepOnWords}, we will let
  $$
  \m_\wd =   \m_{g_1}\ldots\m_{g_n},
  $$ and $$
  \e_\wd = \m_\wd \m_{\wd\inv}.
  $$

  \state Proposition 
  \label RangeVWord
  Given a word $\wd=(g_1,\ldots,g_n)$ in $G$ one has that the range of  
$\m_\wd$, which we will henceforth denote by $\R_\wd$, coincides with
  $$
  \R_\wd = \R_{g_1} \cap \R_{g_1g_2} \cap \ldots\cap\R_{g_1g_2\ldots g_n}.
  $$

  \proof By \lcite{\BasicPrepForWords.ii} we have that
  $$
  \e_\wd =
  \e_{g_1}   \e_{g_1g_2} \ldots \e_{g_1g_2\ldots g_n}.
  $$
  Therefore the range of $\e_\wd$ is precisely
the intersection refered to in the statement.

Since $\e_\wd = \m_\wd\m_{\wd\inv}$ we have that the range of $\e_\wd$
is contained in the range of $\m_\wd$.  By
\lcite{\BasicPrepForWords.i} we have that $\m_\wd = \e_\wd\m_\wd$, so
the reverse inclusion also holds.
  \proofend

\state Proposition
  For every word $\wd$ in $G$ one has that
$(\m_\wd,\m_{\wd\inv})$ is an
interaction.

  \proof
  Clearly $\m_\wd$ and $\m_{\wd\inv}$ are positive bounded linear
maps.
 That
  $$
  \m_\wd\m_{\wd\inv}  \m_\wd = \m_\wd
  \and
  \m_{\wd\inv} \m_\wd\m_{\wd\inv} = \m_{\wd\inv}
  $$
  follows from \lcite{\BasicPrepForWords.i}.
  We must now show that
  $\m_\wd(ab)=\m_\wd(a)\m_\wd(b)$, if either $a$ or $b$ belong to
$\R_{\wd\inv}$.
  Suppose, without loss of generality, that $a\in \R_{\wd\inv}$, and
let $\wd=(g_1,\ldots,g_n)$.
Observing that by \lcite{\RangeVWord} one has that $a\in \R_{g_n\inv}$,
we have
  $$
  \m_\wd(ab) = 
  \m_{g_1}\ldots\m_{g_{n-1}}\m_{g_n}(ab) =
  \m_{g_1}\ldots\m_{g_{n-1}}\Big(\m_{g_n}(a)\m_{g_n}(b)\Big) \$=
  \m_{\wdaux}\Big(\m_{g_n}(a)\m_{g_n}(b)\Big),
  $$
  where $\wdaux=(g_1,\ldots,g_{n-1})$.
  In order to complete the proof we could use induction if we knew
that $\m_{g_n}(a)\in\R_{\wdaux\inv}$.
  In order to prove that this is in fact true observe that by
\lcite{\BasicPrepForWords.ii} we have 
  $$
  \e_{\wdaux\inv}\m_{g_n}(a) =
  \e_{g_{n-1}\inv}  \e_{g_{n-1}\inv g_{n-2}\inv} \ldots
\e_{g_{n-1}\inv g_{n-2}\inv\ldots g_1\inv} \m_{g_n}(a) \$=
  \m_{g_n} \e_{g_n\inv g_{n-1}\inv} \e_{g_n\inv g_{n-1}\inv
g_{n-2}\inv} \ldots \e_{g_n\inv g_{n-1}\inv g_{n-2}\inv\ldots g_1\inv}
(a) \$=
  \m_{g_n} \e_{g_n\inv}
\e_{g_n\inv g_{n-1}\inv} \e_{g_n\inv g_{n-1}\inv
g_{n-2}\inv} \ldots \e_{g_n\inv g_{n-1}\inv g_{n-2}\inv\ldots g_1\inv}
(a) \$=
  \m_{g_n} \e_{\wd\inv} (a) =
  \m_{g_n} (a),
  $$
  thus proving that $\m_{g_n}(a)$ lies in $\R_{\wdaux\inv}$ and hence
completing the proof.
  \proofend

\section{Covariant representations}
  In this section we will consider the natural notion of
representations for interaction groups.

\definition
  \label DefCovariance
  A \stress{covariant representation} of $(A,G,\m)$ in a unital
C*-algebra $B$ is a pair $(\pi,\letrep)$, where $\pi:A\to B$ is a unital
*-homomorphism and $\letrep:G\to B$ is a *-partial representation such
that
  $$
  \rep g \pi(a) \rep{g\inv} = \pi\(\m_g(a)\)   \rep g \rep{g\inv}.
  $$

  Throughout this section we fix a covariant representation
$(\pi,\letrep)$ of the interaction group $(A,G,\m)$ in the C*-algebra
$B$. 

  The following is a generalization of the covariance condition of 
\lcite{\DefCovariance} to words in $G$.

  \state Proposition 
  For every $a\in A$, and for every word $\wd$ in $G$ one has that
  $$
  \rep \wd \pi(a) \rep{\wd\inv} = \pi\big(\m_\wd(a)\big) \rep \wd
\rep{\wd\inv}.
  $$

  \proof
  Let $\wd = (g_1,g_2,\ldots,g_n)$ and put
$\wdaux=(g_1,g_2,\ldots,g_{n-1})$. Then
  $\rep \wd = \rep \wdaux \rep {g_n}$, and hence
  $$
  \rep \wd \pi(a) \rep{\wd\inv} = 
  \rep \wdaux \rep {g_n} \pi(a) \rep {g_n\inv} \rep{\wdaux\inv} =
  \rep \wdaux \pi\(\m_{g_n}(a)\) \idp {g_n} \rep{\wdaux\inv} = \cdots
  $$
  Observing that
  $\idp {g_n} \rep{\wdaux\inv} \={(\BasicPrepForWords.i)}
  \idp {g_n} \idp{\wdaux\inv}\rep{\wdaux\inv} = 
  \idp{\wdaux\inv}  \idp {g_n}  \rep{\wdaux\inv},
  $
  we see that the above equals
  $$
  \cdots =
  \rep \wdaux \pi\(\m_{g_n}(a)\) \rep{\wdaux\inv} \rep{\wdaux}
\idp {g_n} \rep{\wdaux\inv} = \cdots
  $$
  By induction on the length of $\wd$ this equals
  $$
  \cdots =
  \pi\(\m_\wdaux(\m_{g_n}(a))\) \rep{\wdaux} \rep{\wdaux\inv}
\rep{\wdaux} \idp {g_n} \rep{\wdaux\inv} =
  \pi\(\m_\wd(a)\) \rep{\wdaux} \idp {g_n} \rep{\wdaux\inv} \$=
  \pi\(\m_\wd(a)\) \idp {g_1\ldots g_n} \idp{\wdaux}
  \={(\BasicPrepForWords.iii)} 
  \pi\(\m_\wd(a)\) \idp{\wd} =
  \pi\(\m_\wd(a)\) \rep{\wd} \rep{\wd\inv}.
  \proofend
  $$

\state Proposition 
  \label HomosForIdempotnts
  Let $\wd$ be a word in $G$.  Then
  \izitem
  \zitem For every $a$ in $\R_\wd$ one has that $\pi(a)$ commutes
with $\idp \wd$.
  \zitem The map 
  $$
  a\in\R_\wd \mapsto \pi(a)\idp \wd\in B
  $$
  is a (not necessarily unital) *-homomorphism.
  
  \proof
  Let $\wd=(g_1,\ldots,g_n)$.
  In order to prove (i) it is enough to show that $\pi(a)$ commutes with
$\idp {g_1g_2\ldots g_k}$, for all $k\leq n$, given that $\idp
\wd=\idp {g_1}\idp {g_1g_2}\ldots \idp {g_1g_2\ldots g_n}$, by
\lcite{\BasicPrepForWords.ii}.  Noticing that $a\in \R_{g_1g_1\ldots
g_k}$, for all $k$, it suffices to show that for every $h$ in $G$,
and for every $a\in\R_h$, one has that $\pi(a)$ commutes with $\idp h$.
  To see this notice that
  $$
  \idp h \pi(a) \idp h =
  \rep h \rep {h\inv} \pi(a) \rep h \rep {h\inv} \$=
  \pi(\m_h(\m_{h\inv}(a)))  \rep h \rep {h\inv}  \rep h \rep {h\inv}  =
  \pi(a)  \idp h.
  $$
  Applying the same argument to $a^*$ we conclude that $\idp h
\pi(a^*) \idp h = \pi(a^*) \idp h$.  By taking adjoints we then have
that $\idp h \pi(a) \idp h = \idp h \pi(a)$, thus proving (i).
  
Point (ii) then follows immediately from (i).
  \proofend

\state Proposition 
  \label IntroduceMwd
  For each word $\wd$ in $G$ let $\M_\wd$  and $\Kone_\wd$
be the
subsets of $B$ given by 
  $$
  \M_\wd = \overline{\pi(A)\rep \wd \pi(A)}
  \and
  \Kone_\wd = \overline{\pi(A)\idp\wd \pi(A)}.  
  $$
  Then for every $\wd$ one has that
  \izitem 
  \zitem $\M_\wd^* = \M_{\wd\inv}$,
  \zitem $\M_\wd$ is a ternary ring of operators in the sense that
$\M_\wd \M_\wd^* \M_\wd \subseteq \M_\wd$,
  \zitem $\Kone_\wd$ is a closed *-subalgebra of $B$,
  \zitem $\Kone_\wd  \M_\wd \subseteq \M_\wd$,
  \zitem $\M_\wd \Kone_{\wd\inv} \subseteq \M_\wd$,
  \zitem $\overline{\M_\wd \M_\wd^*} = \Kone_\wd$, and
  \zitem $\overline{\M_\wd ^*\M_\wd} = \Kone_{\wd\inv}$.

  \proof We begin by proving (ii).  For this
let $x_1,x_2,x_3\in \M_\wd$ have the simple  form
  $x_i=\pi(a_i)\rep\wd \pi(b_i)$, for i=1,2,3.  Then
  $$
  x_1x_2^*x_3 =
  \pi(a_1)\rep\wd \pi(b_1) \pi(b_2)^*\rep{\wd\inv} \pi(a_2)^*
\pi(a_3)\rep\wd \pi(b_3) \$=
  \pi(a_1)\pi\big(\m_\wd (b_1b_2^*)\big)\rep\wd\rep{\wd\inv}\rep\wd
\pi\big(\m_{\wd\inv}(a_2^*a_3)\big)\pi(b_3) \$=
  \pi(a_1)\pi\big(\m_\wd (b_1b_2^*)\big)\rep\wd \pi\big(\m_{\wd\inv}(a_2^*a_3)\big)\pi(b_3) 
  \in \M_\wd.
  $$
  This proves (ii) which in turn immediately implies (iii) -- (v).  We
leave the easy proof of the other points to the reader.
  \proofend

As a consequence one should notice that $\Kone_\wd$ and $\Kone_{\wd\inv}$ are
Morita--Rieffel equivalent C*-algebras and that $\M_\wd$ is an imprimitivity
bimodule.

\definition
  \label DefineZwd
  Given a word $\wd=(g_1,\ldots,g_n)$ in $G$ we will let
$\Z_\wd$ be the closed linear subspace of $B$ (recall that $B$ is the
codomain of our covariant representation) spanned by the set
  $$
  \{\pi(a_0) \rep{g_1} \pi(a_1) \rep{g_2} \pi(a_2) \ldots \pi(a_{n-1}) \rep
{g_n} \pi(a_n) :
  \forall a_0,\ldots,a_n\in A\}.
  $$
  If $\wd$ is the empty word we set $\Z_\wd=\pi(A)$, by default.

If the reader has indeed understood the subtle difference between the
definitions of $\M_\wd$ and $\Z_\wd$, he or she will recognize that
$\M_\wd\subseteq\Z_\wd$, but not necessarily vice-versa.

  In the next Proposition we will employ  the set $\spr(\wd)$
introduced in \lcite{\DefMuWd}.

  \state Proposition
  \label ZMinM
  If $\wd$ and $\wdaux$ are words in $G$ such that
$\spr(\wd\inv)\subseteq\spr(\wdaux)$ then
  $$
  \Z_\wd\M_\wdaux\subseteq \M_{\wd\wdaux}.
  $$
 
  \def \alt#1{\check#1}   \def \alt#1{#1'}

  \proof 
We need to prove that for every $z\in\Z_\wd$
and $y\in\M_\wdaux$, one has that $zy\in \M_{\wd\wdaux}$. 
  Let $\wd=(g_1,\ldots,g_n)$ and $\wdaux=(h_1,\ldots,h_m)$  and put 
 $\alt\wd=(g_1,\ldots,g_{n-1})$ and $\alt\wdaux=g_n \wdaux$ (concatenation),
so that $\wd\wdaux= \alt\wd\alt\wdaux$.

By density we may clearly suppose that
  $$
  z = \alt z \rep {g_n} \pi(a),
  $$
and
  $$
  y = \pi(b) \rep \wdaux \pi(c),
  $$
  where $\alt z\in\Z_{\alt\wd}$, and $a,b,c\in A$.  Observing that
  $$
  g_n\inv\in \spr(\wd\inv) \subseteq \spr(\wdaux),
  $$
  it follows from \lcite{\BasicPrepForWords.ii} that
$\idp{g_n\inv}\idp\wdaux = \idp\wdaux$, and by
\lcite{\BasicPrepForWords.i}, that $\idp{g_n\inv}\rep\wdaux =
\rep\wdaux$.  So
  $$
  zy = 
  \alt z \rep {g_n} \pi(a) \pi(b) \idp{g_n\inv}\rep\wdaux \pi(c) =
  \alt z \rep {g_n} \pi(ab) \rep{g_n\inv} \rep{g_n}\rep\wdaux
\pi(c) =
  \alt z\pi\big(\m_{g_n}(ab)\big) \rep {g_n} \rep\wdaux \pi(c) \$=
  \alt z \pi\big(\m_{g_n}(ab)\big) \rep{\alt\wdaux} \pi(c) \in
  \Z_{\alt\wd}\M_{\alt\wdaux}.
  \eqno{(\dagger)}
  $$

  On the other hand notice that
  $$
  \spr(\wdaux) \supseteq
  \spr (\wd\inv) = 
  \spr (g_n\inv{\alt\wd}\inv) \={(\SProdConcat)}
  \{1,g_n\inv\} \cup  g_n\inv\spr ({\alt\wd}\inv),
  $$
  so
  $$
  \spr ({\alt\wd}\inv)  \subseteq
  g_n\spr(\wdaux)\ {\buildrel (\SProdConcat) \over\subseteq}\
  \spr(g_n\wdaux) =
  \spr(\alt\wdaux).
  $$
  By induction on $n$ one has that
$\Z_{\alt\wd}\M_{\alt\wdaux}\subseteq \M_{\alt\wd\alt\wdaux} =
\M_{\wd\wdaux}$, so $zy\in\M_{\wd\wdaux}$ by $(\dagger)$.
  \proofend

\state Proposition
  \label IntroduceZX
  Given a finite subset $X\subseteq G$, denote by 
  $$
  \Z^X =\overline{
  \sum_{
  \matrix{\scriptstyle\spr(\wd\inv)\subseteq X \cr \scriptstyle \prd \wd =1}
  }
  \Z_\wd
  }
  $$
  Then $\Z^X$ is a closed *-subalgebra of $B$.

  \proof 
  Let ${\cal W}$ be the set of all words $\wd$ such that
$\spr(\wd\inv)\subseteq X$ and $\prd \wd =1$.  We claim that ${\cal
W}$ is closed both under inversion and concatenation. In fact, given
$\wd\in{\cal W}$ we have by \lcite{\SprOfInverse.ii} that
  $
  \spr((\wd\inv)\inv) = 
  \spr(\wd\inv)   \subseteq X
  $
  and
  $\wprd{\wd\inv} = 1$, so $\wd\inv\in{\cal W}$.  Moreover,
given $\wd_1,\wd_2\in{\cal W}$, observe that
$\wprd{\wd_1\wd_2}= \prd\wd_1 \prd\wd_2 = 1$, and
  $$
  \spr((\wd_1\wd_2)\inv)  =
  \spr(\wd_2\inv\wd_1\inv)  \={(\SProdConcat)} 
  \spr(\wd_2\inv) \cup  \prd{\wd_2\inv} \spr(\wd_1\inv) \$=
  \spr(\wd_2\inv) \cup  \spr(\wd_1\inv)  \subseteq
  X.
  $$
  This proves the claim.  Now, in order to prove the statement it is
enough to realize that
  $$
  \Z_{\wd_1}  \Z_{\wd_2} \subseteq \Z_{\wd_1\wd_2},
  $$
  and
  $$
  (\Z_\wd)^* =   \Z_{\wd\inv}.
  \proofend
  $$

\definition
  \label DefineKwd
  If $\wd$ is a word in $G$ we will let $\K_\wd=\Z^{\spr(\wd)}$.

In other words, $\K_\wd$ is the closed sum of the $\Z_\wdaux$, for all
words $\wdaux$ such that $\prd\wdaux=1$ and $\spr(\wdaux\inv)\subseteq
\spr(\wd)$.

\state Proposition
  \label KwdMwdInMwd
  For every word $\wd$ in $G$ one has that
  $$
  \K_\wd\M_\wd\subseteq\M_\wd.
  $$   

  \proof
  It clearly suffices to prove that $\Z_\wdaux\M_\wd\subseteq\M_\wd$
for all words $\wdaux$ such that $\prd\wdaux=1$ and $\spr(\wdaux\inv)\subseteq
\spr(\wd)$.  Given such a $\wdaux$ observe that
  $$
  \spr(\wdaux) \={(\SprOfInverse.ii)}
  \spr(\wdaux\inv) \subseteq \spr(\wd).
  $$
  Thus we have that $\rep{\wdaux\wd}=\rep\wdaux\rep\wd=\rep\wd$, by
\lcite{\TheFactIMissed}.  Therefore
  $$
  \Z_\wdaux\M_\wd{\buildrel  (\ZMinM) \over\subseteq}
  \M_{\wdaux\wd} =
  \overline{\pi(A)\rep{\wdaux\wd} \pi(A)}  =
  \overline{\pi(A)\rep{\wd} \pi(A)}  =
  \M_\wd.
  \proofend
  $$

  \section{The Toeplitz algebra of an interaction group}
  We shall now  define a ``Toeplitz algebra'' as an
auxiliary step before introducing the main construction, i.e.~the crossed-product algebra.

  \definition
  \label DefToep
  The Toeplitz algebra of an interaction group $(A,G,\m)$ is the
universal unital C*-algebra $\toep$ generated by a copy of $A$ and a
set
  $\{\ts_g: g\in G\}$ subject to the relations
  \izitem
  \zitem $\ts_1=1$,
  \zitem $\ts_g^* = \ts_{g\inv}$,
  \zitem $\ts_g\ts_h\ts_{h\inv} = \ts_{gh}\ts_{h\inv}$,
  \zitem $\ts_g a \ts_{g\inv} = \m_g(a) \ts_g  \ts_{g\inv}$,
  \medskip \noindent for all $a\in A$ and $g,h\in G$.

\medskip
  In order to help the reader to keep track of the objects related to
the Toeplitz algebra we have chosen to decorate these with a hat.
The rule of thumb is therefore that everything that has a hat on it is
related to the Toeplitz algebra.

  Observe that by \lcite{\DefToep.ii-iii} one also has that 
  $\ts_{g\inv}\ts_g\ts_h = \ts_{g\inv}\ts_{gh}$, for all $g$ and $h$ in
$G$.  Therefore, the correspondence 
  $$
  g\in G\mapsto\ts_g\in\toep
  $$
  is a *-partial representation of $G$.
It is therefore clear that:

\state Proposition 
  \label ToeplitzCovRep
  If 
  $$\tj: A \to \toep
  $$
  denotes the canonical map then $(\tj,\ts)$ is a covariant
representation of $(A,G,\m)$, henceforth referred to as the
\stress{Toeplitz covariant representation}.

For future reference we state the universal property of $\toep$ in a
form that emphasizes its role relative to covariant representations:

\state Proposition
  \label PiTimesRep
  If $(\pi,\letrep)$ is a covariant representation of $(A,G,\m)$ in a
C*-algebra $B$ then there exists a unique *-homomorphism $\trxv\pi\letrep:\toep
\to B$ such that
  the following diagram commutes
  $$
  \matrix{
  A \cr\cr
  \tj\, \big\downarrow\ & \searrow\pi \cr\cr
  \toep & \labelarrow{\trxv\pi\letrep} & B \cr\cr
  \ts\, \big\uparrow\ & \nearrow\letrep \cr\cr
  G
  }
  $$

\section{Redundancies}
  \label RedundantSection
  As before we fix a non-degenerated interaction group $(A,G,\m)$.  We
will now make extensive use of the spaces $\M_\wd$ of
\lcite{\IntroduceMwd}, the $\Z_\wd$ of \lcite{\DefineZwd}, and the
$\K_\wd$ of \lcite{\DefineKwd}, all of these relative to the Toeplitz
covariant representation.  According to our rule of thumb we will
denote these by $\tM_\wd$, $\tZ_\wd$, and $\tK_\wd$, respectively.


\definition 
  \label DefineRedundancies
  Let $\wd$ be a word in $G$.
By an \stress{$\wd$-redundancy} we shall mean any element $k\in\tK_\wd$
such that $k\tM_\wd=\{0\}$.  This is of course equivalent to
  $$
  k\,\tj(b)\,\ts_\wd=0
  \for b\in A.
  $$

Observe that, unlike
  \scite{\larsen}{1.2} and
  \scite{\fowler}{2.5},
  our notion of redundancies involves many group elements at once in
the sense that the definition of $\tK_\wd$ involves many $\ts_g$.

The following is the main concept introduced by this article:

\definition 
  \label DefineCp
  Let $(A,G,\m)$ be an interaction group.  The crossed-product of $A$
by $G$ under $\m$ is the C*-algebra $A\crossproduct_\m G$,
obtained by taking the quotient of\/ $\toep$ by the ideal generated by
all redundancies.  The quotient map will be denoted
  $$
  q : \toep \to \cp,
  $$
  the composition
  $$
  A \labelarrow \tj \toep \labelarrow q \cp
  $$
  will be denoted by $\cpj$,
  and for every $g$ in $G$ we will let
  $$
  \cps_g = q(\ts_g).
  $$
  It is therefore evident that $(\cpj,\cps)$ is a covariant
representation of $(A,G,\m)$ in $\cp$, henceforth referred to as the
\stress{fundamental covariant representation}.

Recall from \lcite{\PiTimesRep} that every covariant representation
$(\pi,\letrep)$ induces a *-homomorphism $\trxv\pi\letrep$ defined
on $\toep$.  There is no reason, however, for  $\trxv\pi\letrep$
to vanish on redundancies.

\definition
  \label IntroduceStrongly
  Let $(\pi,\letrep)$ be a covariant representation of $(A,G,\m)$ in
a C*-algebra $B$.
  We will say that $(\pi,\letrep)$ is
\stress{strongly covariant} if $\trxv\pi\letrep$ vanishes on all
redundancies.  In this case one therefore has a *-homomorphism
$\rxv\pi\letrep:\cp\to B$ such that the following diagram commutes
  $$
  \matrix{
  A \cr\cr
  \cpj\, \big\downarrow\ & \searrow\pi \cr\cr
  \cp & \labelarrow{\rxv\pi\letrep} & B \cr\cr
  \cps\, \big\uparrow\ & \nearrow\letrep \cr\cr
  G
  }
  $$

It is evident that the fundamental covariant representation
$(\cpj,\cps)$ is strongly covariant.

\section{Grading of the crossed product}
  Recall that a C*-algebra $B$ is said to be graded over the group $G$
if one is given a independent family $\{B_g\}_{g\in G}$ of closed
linear subspaces of $B$, called the \stress{grading subspaces}, such that
$\bigoplus_{g\in G} B_g$ is dense in $B$ and such that for every $g$
and $h$ in $G$ one has that $B_g B_h \subseteq B_{gh}$ and
$B_g^*=B_{g\inv}$.  
  According to \scite{\amena}{3.4} a grading is said to be a
\stress{topological grading} if there exists a conditional expectation
from $B$ onto $B_1$ which vanishes on $B_g$, for all $g\neq1$.

In this section we will prove that $\cp$ admits a canonical
topological grading
such that for all $g$ in $G$, one has that $\cps_g$ lies in the
grading subspace associated to $g$.

In the following result we will employ the  full C*-algebra of $G$,
denoted $C^*(G)$, as well as the universal representation of $G$
  $$
  u: G \to C^*(G).
  $$

\state Proposition
  \label Amplification
  Let $(\pi,\letrep)$ be a covariant representation of $(A,G,\m)$ in
a C*-algebra $B$.  Consider the maps
  $$
  \pi'=\pi\* id : A \to B\* C^*(G)
  $$
  and
  $$
  \letrep' = \letrep\*u: G \to  B\* C^*(G).
  $$
  Then $(\pi',\letrep')$ is a covariant representation.  Moreover if
$(\pi,\letrep)$ is strongly covariant then so is $(\pi',\letrep')$.
We will refer to $(\pi',\letrep')$ as the \stress{amplification} of
$(\pi,\letrep)$.

  \proof
  We leave everything up to the reader, except for the last statement.
If $\wdaux=(h_1,\ldots,h_n)$ is a word in $G$ and
  $z$
  is any element in $\tZ_\wdaux$
  notice that 
  $$
  \trxv{\pi'}{\letrep'}(z) =
  \trxv{\pi}{\letrep}(z)\*u_h,
  $$ 
  where $h=\prd\wdaux$.  In case $\prd\wdaux=1$ we therefore have that 
  $$
  \trxv{\pi'}{\letrep'} = \trxv{\pi}{\letrep}\*1
  $$
  on $\tZ_\wdaux$.  This will clearly also be the case on $\tK_\wd$ for
every word $\wd$.  Given any  $\wd$-redundancy $k\in \tK_\wd$  we
conclude that
  $$
  \trxv{\pi'}{\letrep'}(k) = (\trxv{\pi}{\letrep})(k)\*1.
  $$
  If $(\pi,\letrep)$ is strongly covariant then 
  $\trxv{\pi}{\letrep}(k)=0$, for all redundancies $k$ and hence the
same is true for $\trxv{\pi'}{\letrep'}(k)$.
  \proofend

\state Proposition
  \label CpGrading
  $\cp$ admits a topological $G$-grading $\{C_g\}_{g\in G}$
such that $\cpj(A)\subseteq C_1$ and $\cps_g\in C_g$, for every $g$ in
$G$.

  \proof  
  Throughout this proof whenever we speak of the $\Z_\wd$ of
\lcite{\DefineZwd} it will be with respect to the fundamental
covariant representation (see \lcite{\DefineCp}).
  For every $g$ in $G$ let
  $$
  C_g = \overline{
  \sum_{\prd\wd=g}
  \Z_\wd}.
  $$
  We leave it for the reader to prove that $\{C_g\}_{g\in G}$
satisfies \scite{\amena}{3.3.i--iii}.  Next we will provide a
bounded linear map
  $$
  F:\cp\to C_1
  $$
  which is the identity on $C_1$ and which vanishes on all $C_g$ for $g\neq1$.

For this consider the fundamental covariant representation $(\cpj,\cps)$ and
let $(\cpj',\cps')$ be the corresponding amplification given by
\lcite{\Amplification}.  Since $(\cpj,\cps)$ is strongly covariant
then so is $(\cpj',\cps')$, and hence by \lcite{\IntroduceStrongly}
we get
  $$
  \rxv{\cpj'}{\cps'} : \cp \to (\cp)\*C^*(G).
  $$
  It is elementary to verify that 
  $$
  \rxv{\cpj'}{\cps'}(z) =  z\*u_g
  \for z\in C_g.
  $$

Let $tr$ be the standard trace on $C^*(G)$ and set 
  $$
  F=  (id\*tr)\circ (\rxv{\cpj'}{\cps'}).
  $$
  For $z\in C_g$ we therefore have
  $$
  F(z) = (id\*tr)(z\*u_g) = \left\{\matrix {z, & \hbox{ if }g=1, \cr\cr 0, &
\hbox{ if }g\neq1.}\right.
  $$
  We therefore see that $F$ satisfies the required properties and
hence we may invoke \scite{\amena}{3.3} to conclude that $\{C_g\}_{g}$
is a topological grading for $\cp$.

  Finally it is easy to see that $\cpj(A)\subseteq C_1$, and
$\cps_g\in\Z_g\subseteq C_g$, as desired.
  \proofend

  \section{Several Hilbert modules}
  \label HilbertModSection
  So far nothing guarantees that $\toep$ or $\cp$ contain a single
non-zero element!  In order to show that these algebras have any
substance at all we need to do some work.  This section is dedicated
to obtaining several technical results in preparation for the
description of certain representations of $\toep$ and $\cp$.

  \state  Proposition
  \label Inequal
  For all $a$ in $A$ one has that $\m_g(a^*)\m_g(a) \leq \m_g(a^*a)$.

  \proof We have
  $$
  \m_g(a^*)\m_g(a) =
  \m_g\(\e_{g\inv}(a^*)\)  \m_g\(\e_{g\inv}(a)\) \$=
  \m_g\(\e_{g\inv}(a^*)\e_{g\inv}(a) \)\leq
  \m_g\(\e_{g\inv}(a^*a)\)=
  \m_g(a^*a),
  $$
  where the inequality follows from \scite{\takesaki}{IV.3.4}.
  \proofend

Recall that by \lcite{\BasicPreps.iii} the $\e_g$ commute amongst
themselves.  It then easily follows that the composition of any number
of $\e_g$'s is again a conditional expectation onto the intersection
of the $\R_g$'s involved.

  \definition
  \label ExRx
  For each finite subset $X\subseteq G$ we will let 
  $$
  \R_X=\bigcap_{g\in X}\R_g
  \and
  \e_X=\prod_{g\in X}\e_g,
  $$
  so that $\e_X$ is a conditional expectation onto $\R_X$. 

Since each $\e_g$ is non-degenerated we clearly have that the $\e_X$
are non-degenerated as well. 

\state Proposition  For every $g$ in $G$ and every finite subset
$X\subseteq G$ we have that
  \izitem
  \zitem $\m_g(\R_X)\subseteq \R_{gX}$,
  \zitem if\/ $1\in X$ then the inclusion in (i) becomes an equality.

  \proof
  Observing that 
  $$
  \m_g\e_h = \e_{gh}\m_g,
  $$
  by \lcite{\BasicPreps.ii}  we have
  $$
  \m_g(\R_X) =
  \m_g\(\prod_{h\in X}\e_h(A)\) =
  \prod_{h\in X}\e_{gh}(\m_g(A)) \subseteq
  \prod_{h\in X}\e_{gh}(A) = \R_{gX}.
  $$
  If $1\in X$ then
  $$
  \R_{gX} =
  \prod_{h\in X}\e_{gh}(A) =
  \prod_{h\in X}\e_{gh}\e_g(A) =
  \prod_{h\in X}\e_{gh}\m_g\m_{g\inv}(A) \$\subseteq
  \prod_{h\in X}\e_{gh}\m_g(A) =
  \m_g\(\prod_{h\in X}\e_h(A)\) =
  \m_g\(\R_X\).
  \proofend
  $$

  \definition For every finite subset $X\subseteq G$ we will let 
$\H_X$ be the right Hilbert $\R_X$--module obtained by
completing $A$ under the $\R_X$--valued inner-product defined by
  $$
  \<a,b\>_X = \e_X(a^*b)
  \for a,b\in A.
  $$

Notice that the standing hypothesis according to which the $\e_g$ are
non-degenerated implies that the above inner-product is
likewise non-degenerated.

  \state Proposition
  \label AdjointMap
 For every $g$ in $G$, and every finite subset
$X\subseteq G$, there exists a bounded linear map
  $$
  \n_g : \H_X \to \H_{gX}
  $$
  such that 
  $
  \n_g(a) = \m_g(a),
  $
  for all $a$ in $A$ (notice that we are not using any special
decoration to denote the image of an element $a\in A$ within  $\H_X$).
  If moreover $1\in X$ then
  $$
  \<\n_g(\xi),\eta\>_{gX}=
  \m_g\(\<\xi,\n_{g\inv}(\eta)\>_X\)
  \for \xi\in \H_X \for\eta\in \H_{gX}.
  $$

  \proof
  For $a\in A$ we have
  $$
  \<\m_g(a),\m_g(a)\>_{gX} =
  \e_{gX}\(\m_g(a^*)\m_g(a)\) {\buildrel(\Inequal) \over \leq}
  \e_{gX}\(\m_g(a^*a)\) =
  \m_g\(\e_X(a^*a)\) =
  \m_g\(\<a,a\>_X\),
  $$
  so that $\|\m_g(a)\|_{gX} \leq \|a\|_X$ and hence the correspondence
$a \mapsto \m_g(a)$ extends to  a bounded linear map $\n_g:\H_X\to  \H_{gX}$
such that $\n_g(\inc a) = \inc {\m_g(a)}$ , for all $a$ in $A$.

Assuming that $1\in X$ we have that $g\in gX$ and hence $\e_{gX} =
\e_{gX} \e_g$, so that given $a,b\in A$ we have
  $$
  \<\m_g(a),b\>_{gX} = 
  \e_{gX}\(\m_g(a^*)b\) =
  \e_{gX}\(\e_g\(\m_g(a^*)b\)\) \$=
  \e_{gX}\(\m_g(a^*)\e_g(b)\) =
  \e_{gX}\(\m_g(a^*)\m_g(\m_{g\inv}(b))\) \$=
  \e_{gX}\(\m_g(a^*\m_{g\inv}(b))\) =
  \m_g\(\e_X(a^*\m_{g\inv}(b))\) =
  \m_g\(\<a,\m_{g\inv}(b)\>_X\).
  $$
  The conclusion then follows from the density of $\inc A$ both in
$\H_X$ and in $\H_{gX}$.
  \proofend

  \definition
  \label DefLing
  Given $g$ in $G$, let $\LL gX$ denote the set of all maps $T:\H_X\to
\H_{gX}$, for which there exists a map $S:\H_{gX}\to \H_X$ satisfying
  $$
  \<T\xi,\eta\>_{gX}=
  \m_g\(\<\xi,S\eta\>_X\)
  \for \xi\in \H_X \for\eta\in \H_{gX}.
  $$

Observe that because $\m_g$ preserves adjoint, the equation displayed above is equivalent to 
  $$
  \<\eta,T\xi\>_{gX}=
  \m_g\(\<S\eta,\xi\>_X\)
  \for \xi\in \H_X \for\eta\in \H_{gX}.
  $$
  Observe also that when $1\in X$, the map $\n_g$ of
\lcite{\AdjointMap} lies in $\LL gX$.

  \state Proposition
  \label IntroduceTstar
  When $T\in \LL gX$ and $g\inv\in X$, one has that the map $S$
mentioned in \lcite{\DefLing} is uniquely determined.
We therefore denote
it by $T^*$ and call it  the \stress{adjoint of $T$}.
 
  \proof
  We claim that if  $\zeta\in \H_X$ is such that
  $$
  \m_g\(\<\xi,\zeta\>_X\) = 0
  \for \xi\in \H_X,
  $$
  then $\zeta=0$. In order to prove this claim plug 
$\xi=\zeta$ above so that
  $$
  0 =
  \m_{g\inv}\m_g\(\<\zeta,\zeta\>_X\) =
  \e_{g\inv}\(\<\zeta,\zeta\>_X\) =
  \<\zeta,\zeta\>_X,
  $$
  where the last equality follows from the fact that the range of
$\<\cdot,\cdot\>_X$ is contained in $\R_X$ which in turn is contained
in $\R_{g\inv}$ because $g\inv\in X$.  Therefore $\zeta=0$.

If $S_1$ and $S_2$ are maps satisfying the conditions of
\lcite{\DefLing} we have for all $\xi\in \H_X$ and $\eta\in \H_{gX}$
that
  $$
  0 = \m_g\(\<\xi,S_1(\eta)-S_2(\eta)\>_X\),
  $$
  so that $S_1(\eta)=S_2(\eta)$ by the claim.
  \proofend

Combining \lcite{\AdjointMap} and \lcite{\IntroduceTstar} and suposing
that $\{1,g\inv\}\subseteq X$ we have that
$(\n_g)^*=\n_{g\inv}$. 

  \state Proposition  
  \label StarSaysAll
  Any $T\in \LL gX$ is bounded and $\C$-linear.  If moreover $g\inv\in X$
then 
  $$
  T(\xi a) = T(\xi)\m_g(a)
  \for\xi\in \H_X \for a\in \R_X.
  $$

  \proof
  Given $\lambda\in\C$, $\xi_1,\xi_2\in \H_X$, and $\eta\in \H_{gX}$ we have
  $$
  \<T(\xi_1+\lambda\xi_2),\eta\>_{gX}=
  \m_g\(\<\xi_1+\lambda\xi_2,T^*\eta\>_X\) \$=
  \m_g\(\<\xi_1,T^*\eta\>_X\) + \bar \lambda\m_g\(\<\xi_2,T^*\eta\>_X\) =
  \<T(\xi_1),\eta\>_{gX} +\bar\lambda \<T(\xi_2),\eta\>_{gX}\$=
  \<T(\xi_1)+\lambda T(\xi_2),\eta\>_{gX},
  $$
  whence 
  $T(\xi_1+\lambda\xi_2)= T(\xi_1)+\lambda T(\xi_2)$.

  In order to prove that $T$ is bounded suppose that $\lim_n\xi_n= 0$,
while $\lim_nT(\xi_n) = \zeta$.  Then
  $$
  \<\zeta,\zeta\>_{gX} =
  \lim_n \<T(\xi_n),\zeta\>_{gX} =
  \lim_n \m_g\(\<\xi_n,T^*\zeta\>_X\) =0,
  $$
  which implies that $\zeta=0$, and hence the conclusion follows by the Closed
Graph Theorem.

  Given $\xi\in \H_X$, $\eta\in \H_{gX}$, and $a\in \R_X$ we have that
  $$
  \<T(\xi a),\eta\>_{gX}=
  \m_g\(\<\xi a,T^*\eta\>_X\) =
  \m_g\(a^*\<\xi,T^*\eta\>_X\) = \cdots
  $$
  Assuming that $g\inv\in X$ we have that
$\<\xi,T^*\eta\>_X\in \R_{g\inv}$ and so the above equals 
  $$
  \cdots =
  \m_g(a^*)\m_g\(\<\xi,T^*\eta\>_X\) =
  \m_g(a^*)\<T\xi,\eta\>_{gX} =
  \<T(\xi) \m_g(a),\eta\>_{gX},
  $$
  thus proving that 
  $T(\xi a)=T(\xi) \m_g(a)$.  
  \proofend

  \state Proposition
  \label BasicStar
  Suppose that $\{1,g\inv\}\subseteq X$.   Then for every $T\in\LL gX$ one
has that 
  \izitem
  \zitem $T^*\in\LL{g\inv}{gX}$,
  \zitem $(T^*)^*=T$,
  \zitem $\|T^*\|=\|T\|$,
  \zitem $\|T^*T\|=\|T\|^2$.

  \proof Initially observe that the hypothesis that $1\in X$ says that
$g\in gX$, in which case the adjoint of an operator in $\LL{g\inv}{gX}$
is well defined by \lcite{\IntroduceTstar}.

 Given that $g\inv\in X$ we have that $\e_X = \e_{g\inv} \e_X$
and hence $\R_X\subseteq \R_{g\inv}$.  Therefore for any 
  $\xi\in \H_X$ and 
  $\eta\in \H_{gX}$ 
  we have that $\<\xi,T^*\eta\>_X
\in \R_{g\inv}$ so
  $$
  \<\xi,T^*\eta\>_X = 
  \m_{g\inv}\m_g\(\<\xi,T^*\eta\>_X\) = 
  \m_{g\inv}\(\<T\xi,\eta\>_{gX}\).
  $$
  This proves (i) and (ii).  Next observe that  for all $\xi\in \H_X$
we have
  $$
  \|T\xi\|^2 =
  \|\<T\xi,T\xi\>_{gX}\|  = 
  \|\m_g\(\<T^*T\xi,\xi\>_X\)\| \leq
  \|T^*T\|\,\|\xi\|^2,
  $$
  proving that $\|T\|^2\leq \|T^*T\|$.  From this  (iii)
and (iv) follow without difficulty.
  \proofend

  \state Proposition 
  \label Composition
  Supose that 
    $\{g\inv,g\inv h\inv\}\subseteq X$ 
  and let $T\in\LL gX$ and $S\in\LL h{gX}$.  Then $ST\in\LL{hg}X$ and
$(ST)^*=T^*S^*$.

  \proof  
  Initially observe that $h\inv\in gX$ whence $S^*$ is well defined by
\lcite{\IntroduceTstar}.   
For $\xi\in \H_X$ and $\zeta\in \H_{hgX}$ we have
  $$
  \<ST\xi,\zeta\>_{hgX} =
  \m_h\(\<T\xi,S^*\zeta\>_{gX}\) =
  \m_h\(\m_g\(\<\xi,T^*S^*\zeta\>_{X}\)\) = \cdots
  $$
  Letting $a=\<\xi,T^*S^*\zeta\>_{X}$ we have that $a\in \R_X$ and
hence $a = E_X(a) = E_{g\inv}E_X(a) \in \R_{g\inv}$ so the above
equals
  $$
  \cdots =
  \m_h\m_g(a)  =
  \m_h\m_g\m_{g\inv}\m_g(a)  =
  \m_{hg}\m_{g\inv}\m_g(a) =
  \m_{hg}(a)  =
  \m_{hg}\(\<\xi,T^*S^*\zeta\>_{X}\),
  $$
  proving the statement.
  \proofend
  
  \state Proposition
  \label TTsPositive
  If\/ $\{1,g\inv\}\subseteq X$ then for every $T\in \LL gX$ one has that $T^*T$ is
a positive element in the C*-algebra $\Lin(\H_X)=\LL 1X$.

  \proof
  First note that since $g\inv\in X$ we have by \lcite{\BasicStar.i}
that $T^*\in \LL{g\inv}{gX}$.  Applying \lcite{\Composition} (with $h=g\inv$) we
then conclude that $T^*T\in\LL 1X$.
  For $\xi\in \H_X$ we have
  $$
  \<T^*T\xi,\xi\>_X =
  \m_{g\inv}\(\<T\xi,T\xi\>_{gX}\) \geq 0,
  $$
  and hence $T^*T$ is positive.
  \proofend


  Let $X$ and $Y$ be finite subsets of $G$ such that $X\subseteq Y$.
Then $E_Y = E_Y E_X$, so for every $a\in A$ we have
  $$
  \|\<a,a\>_Y\| =
  \|E_Y(a^*a)\| =
  \|E_YE_X(a^*a)\| \leq
  \|E_X(a^*a)\| =
  \|\<a,a\>_X\|,
  $$
  and hence the correspondence $a\in A \mapsto \inc a\in \H_Y$ extends
to give a contractive\fn{We use the term ``contractive'' to mean
``non-expansive'', i.e.~that $\|\iota(\xi)\|\leq \|\xi\|$, for all
$\xi$ in $\H_X$.} linear map
  $$ 
  \iota : \H_X \to \H_Y,
  $$
  whose restriction to $A$ (or rather to the canonical dense copy of
$A$ within $\H_X$) is the identity.

Because
  $$
  \<\iota (a),\iota(a)\>_Y =
  E_Y(a^*a) =
  E_YE_X(a^*a) =
  E_Y\(\<a,a\>_X\)
  \for a\in A,
  $$
  we have that 
  $$
  \<\iota(\xi),\iota(\xi)\>_Y =
  E_Y\(\<\xi,\xi\>_X\)
  \for \xi\in \H_X.
  \eqno {(\seqnumbering)}
  \label RelationInnerProd
  $$

Since $\e_Y$ is non-degenerated one sees that $\iota$ is an injective
map.  It will often be convenient to think of $\H_X$ as a dense
subspace of $\H_Y$.  Nevertheless care must be taken to account for
the fact that $\H_X$ and $\H_Y$ are quite different objects, having
different norms and being Hilbert modules over different algebras.

  \state Proposition 
  \label Extension
  Let $X$ and $Y$ be finite subsets of $G$, let $g\in G$ be such
that $\{1,g\inv\}\subseteq X\subseteq Y$, and let $T\in\LL gX$.  Then
there exists a unique bounded operator $\tilde T: \H_Y\to \H_{gY}$ such
that the diagram
  $$
  \matrix{
    \H_X & \labelarrow{T} & \H_{gX}\cr
    \vrule height 13pt width 0pt \iota \downarrow\ && \ 
\downarrow\iota\cr
    \H_Y & \labelarrow{\tilde T} & \H_{gY}}
  $$
  commutes.  Moreover $\tilde T\in\LL gY$, $\|\tilde T\|\leq\|T\|$,
and $(\,\tilde T\,)^* = \tilde{T^*}$.

  \proof 
  We should first observe that we are denoting by $\iota$ both
inclusions $\H_X\to \H_Y$ and $\H_{gX}\to \H_{gY}$, leaving for the context
to distinguish which is which.

  By 
  \lcite{\TTsPositive}
  we have that $T^*T\in \LL 1X$ and by \lcite{\BasicStar.iv} it
follows that $T^*T\leq \|T\|^2$.  Therefore for all $\xi\in \H_X$ we
have
  $$
  \<T\xi,T\xi\>_{gX} = 
  \m_g\(\<T^*T\xi,\xi\>_X\) \leq
  \|T^2\|\,\m_g\(\<\xi,\xi\>_X\).
  $$
  So by \lcite{\RelationInnerProd} we have that
  $$
  \<\iota T\xi,\iota T\xi\>_{gY} = 
  E_{gY}\(\<T\xi,T\xi\>_{gX}\) \leq
  \|T^2\|\,E_{gY}\m_g\(\<\xi,\xi\>_X\) \$=
  \|T^2\|\,\m_g E_Y\(\<\xi,\xi\>_X\) =
  \|T^2\|\,\m_g\(\<\iota \xi,\iota \xi\>_Y\),
  $$
  from where one deduces that 
  $
  \|\iota T\xi\|\leq \|T\|\,\|\iota \xi\|.
  $
  Given that $\iota(\H_X)$ is dense in $\H_Y$ we have that 
  the correspondence $\iota\xi \mapsto \iota T\xi$  extends
to a bounded linear map $\tilde T: \H_Y \to \H_{gY}$ such that  the diagram
above commutes and such that  $\|\tilde T\|\leq\|T\|$.  Again because
$\iota(\H_X)$ is dense in $\H_Y$ we have that $\tilde T$ is  uniquely determined.

Since $\{1,g\}\subseteq gX$, the above reasoning applies to $T^*$ so we
may speak of $\tilde{T^*}$ as well.
  Given $\xi\in \H_X$ and $\eta\in \H_{gX}$
we have
  $$
  \<\tilde T\iota\xi,\iota \eta\>_{gY} =
  \< \iota T\xi,\iota\eta\>_{gY} =
  \e_{gY}\(\< T\xi,\eta\>_{gX}\) =
  \e_{gY}\(\m_g\(\<\xi,T^*\eta\>_X\)\) \$=
  \m_g\(\e_Y\(\<\xi,T^*\eta\>_X\)\) =
  \m_g\(\<\iota\xi,\iota T^*\eta\>_Y\)=
  \m_g\(\<\iota\xi,\tilde{T^*}\iota\eta\>_Y\).
  $$
  Since  $\iota(\H_X)$ is dense in $\H_Y$ and $\iota(\H_{gX})$ is dense in $\H_{gY}$ it
follows that 
  $$
  \<\tilde T\xi,\eta\>_{gY} =
  \m_g\(\<\xi,\tilde{T^*}\eta\>_Y\)
  \for \xi\in \H_Y \for \eta \in \H_{gY},
  $$
  whence $\tilde T\in\LL gY$ and
$(\tilde T)^* = \tilde{T^*}$.
  \proofend

  \state Proposition 
  \label ExtendCompose
  Suppose that $\{1,g\inv,g\inv h\inv\}\subseteq X\subseteq Y$.  Then
for every 
$T\in\LL gX$ and $S\in\LL h{gX}$ one has that  $\tilde{(ST)} = \tilde S
\tilde T$.

  \proof
  Initially observe that, by abuse of language, we
are denoting by tilde the correspondences given by \lcite{\Extension}
in the all three cases:

\bigskip  $\LL gX\,\ \to \ \,\LL gY$ (for $T$),

\medskip  $\LL h{gX} \to \LL h{gY}$ (for $S$), and

\medskip  $\LL{hg}X \to \LL{hg}Y$ (for $ST$).
  
\bigskip\noindent Also notice that $ST\in\LL{hg}X$ by
\lcite{\Composition}.
  Considering the diagram
  $$
  \matrix{
    \H_X & \labelarrow{T} & \H_{gX} & \labelarrow{S} & \H_{hgX}\cr
    \vrule height 13pt width 0pt \downarrow && \downarrow && \downarrow \cr
    \H_Y & \labelarrow{\tilde T} & \H_{gY} & \labelarrow{\tilde S} & \H_{hgY}}
  $$
  we deduce that
  $\tilde{(ST)}=\tilde S \tilde T$ by the uniqueness part of
\lcite{\Extension}.
  \proofend

  If we specialize \lcite{\Extension} and \lcite{\ExtendCompose} to
the case in which $g=h=1$, and if we suppose that $1\in X\subseteq Y$,
we conclude that the correspondence
  $$
  T\in \LL 1X \mapsto \tilde T \in \LL 1Y
  $$
  is a C*-algebra homomorphism.

Moreover, since $\tilde T$ is an extension of $T$, it is clear that
this homomorphism is injective and hence isometric.

  \state Corollary
  \label IsometricExt
  If\/ $\{1,g\inv\}\subseteq X\subseteq Y$ and $T\in\LL gX$ then
  $\|\tilde T\| = \|T\|$.

  \proof
  By \lcite{\TTsPositive} it follows that $T^*T\in\Lin(\H_X)$.  Since
$1\in X$ we have by the discussion above that
  $$
  \|T^*T\| = \|\widetilde{(T^*T)}\| = \|\tilde T ^* \tilde T\|.
  $$
  The conclusion then follows from \lcite{\BasicStar.iv}.
  \proofend

\section{The regular covariant representation}
  In this section we will describe representations of $\toep$ and
$\cp$ which will, among other things, show that the maps $\tj$ and
$\cpj$ are injective.

  Given $g$ in $G$, let $\I_g$ be the subset of the power set of $G$
defined by 
  $$ 
  \I_g  = \big\{X\subseteq G: 1,g\inv\in X\big\}.
  $$
  Ordered by inclusion we have that $\I_g$ is a directed set.
For each $X,Y\in\I_g$ let us denote by $\ext YX$ the correspondence
  $$
  T\in \LL gX \to \tilde T \in \LL gY
  $$
  of \lcite{\Extension}, which we now know is an isometry.

  We want to view the collection of Banach
spaces $\{\LL gX\}_{X\in\I_g}$, equipped with the  maps
$\{\ext YX\}_{X,Y\in\I_g}$, as a directed system.
  In order to justify this we must show that, whenever $X\subseteq
Y\subseteq Z$, one has that $\ext ZY \ext YX = \ext ZX$.
  But this is easy to see once one realizes that the diagram
  $$
  \matrix{
    \H_X & \labelarrow{T} & \H_{gX}\cr
    \vrule height 13pt width 0pt \downarrow && \downarrow\cr
    \H_Y & \labelarrow{\tilde T} & \H_{gY} \cr
    \vrule height 13pt width 0pt \downarrow && \downarrow\cr
    \H_Z & \labelarrow{\tilde {\tilde T}} & \H_{gZ}
  }
  $$
  commutes, where $\tilde T = \ext YX(T)$, and $\tilde {\tilde T} =
\ext ZY(\tilde T)$.

In the next Definition we shall refer to the inductive limit in the
category of Banach spaces with isometric linear maps as morphisms.  It
is easy to prove that this category does indeed admit inductive
limits.

  \definition Given $g\in G$ we shall denote the inductive limit of
the directed family
  $\big\{\LL gX\big\}_{X\in\I_g}$ by $B_g$.  
  The dense subspace of $B_g$ formed by the union of the $\LL gX$, as
$X$ ranges in $\I_g$, will be denoted by $B^0_g$.
  
When $X=\{e\}$ it is clear that $\H_X=A$, and $\LL 1X$ is *-isomorphic
to $A$.  We therefore have an injective *-homomorphism
  $$
  i:A\simeq\LL 1X \hookrightarrow B_1,
  \eqno {(\seqnumbering)}   
  \label AinBOne
  $$
  which we will often use to  view $A$ as a subalgebra of $B_1$.  

Given $T\in B_g^0$ and $S\in B_h^0$ choose $X\in \I_g$ and $Y\in
\I_h$ so that $T\in \LL gX$ and $S\in \LL hY$.
We may assume without loss of generality that $Y=gX$ (if not one may 
replace $X$ by $X\cup g\inv Y$ and $Y$ by $Y\cup g X$).
By \lcite{\Composition} we have that $ST\in\LL{hg}X\subseteq B_{hg}$.

Suppose that instead of choosing the above $X$ and $Y$ we take
$X'\supseteq X$ and $Y'\supseteq Y$, still satisfying the relation
$Y'=gX'$.  In the inductive limit one would then identify $T$ with its
extension $\tilde T\in \LL g{X'}$, and $S$ with $\tilde S\in \LL
h{Y'}$.  By \lcite{\ExtendCompose} it follows that $\widetilde{(ST)} =
\tilde S \tilde T$ which means that the class of $ST$ in $B_{hg}$ does
not depend on the above choices.

We have therefore defined an operation
  $$
  B_h^0 \times B_g^0 \to B_{hg}^0
  $$
  which obviously satisfies $\|ST\|\leq \|S\|\,\|T\|$ and hence
extends to a \stress{multiplication  operation}
  $$
  B_h \times B_g \to B_{hg}.
  $$
  It is also clear that the adjoint operation introduced in
\lcite{\IntroduceTstar} gives an involution
  $$
  T\in B_g \mapsto T^*\in B_{g\inv}.
  $$

  \state Proposition
  The collection $\B = \big\{B_g\}_{g\in G}$ is a Fell
bundle \cite{\fell} with the operations defined above.

  \proof
  Follows from the results of section \lcite{\HilbertModSection}.
  \proofend

Given $g\in G$, observe that by \lcite{\AdjointMap} we have that
$\n_g\in\LL g{X}$, for any  $X\in\I_g$.  Henceforth we
will identify $\n_g$ with its image in $B_g$ and hence we will think
of $\n_g$ as an element of $\CB$, the cross sectional C*-algebra of
$\B$ \cite{\fell}.

  \state Proposition
  \label TheMainPrep
  The map 
  $$
  g\in G\mapsto \n_g\in  \CB
  $$
  is a partial representation of\/ $G$ in $\CB$ and  for every $a\in
A\subseteq B_1 \subseteq \CB$
we have that
  $$
  \n_g a \n_{g\inv} = \m_g(a) \n_g  \n_{g\inv}
  \for g\in G.
  $$

  \proof
  Let $X=\{1,g\}$ and let us view $\n_{g\inv}\in\LL {g\inv}{X}$,
$\n_g\in\LL g{g\inv X}$, $a\in \LL 1{g\inv X}$, and $\m_g(a)\in \LL
1{X}$.  Then, proving the displayed relation above amounts to
comparing the following compositions
  $$
  \H_X \labelarrow{\n_{g\inv}}  
  \H_{g\inv X} \labelarrow{\n_g} 
  \H_X \labelarrow{\m_g(a)}  \H_X
  $$
  and
  $$
  \H_X \labelarrow{\n_{g\inv}}  
  \H_{g\inv X} \labelarrow{a}
  \H_{g\inv X} \labelarrow{\n_g} \H_X.
  $$
  Since $A$ is dense in $\H_X$ it is enough to prove that these compositions
agree on $A$.  Given any  $b\in A$ we have
  $$
  \n_g a \n_{g\inv} (b) = 
  \m_g (a \m_{g\inv} (b)) =
  \m_g (a) \m_g(\m_{g\inv} (b)) =
  \m_g(a) \n_g  \n_{g\inv}(b).
  $$
  That $\n$  is a partial representation follows from a similar
  argument.
  \proofend

An equivalent way to state the above result is:

\state Corollary
  \label BundleCovarRep
  Let \ $i:A\to B_1\subseteq C^*(\B)$ \ be the canonical inclusion
described in \lcite{\AinBOne}.  Then $(i,\n)$ is a covariant
representation of $(A,G,\m)$ in $C^*(\B)$.  This will be called the
\stress{regular covariant representation}.

Another important  consequence is:

  \state Corollary
  \label AEmbedsInToep
  The canonic map \ $\tj:A\to \toep$ is injective.

  \proof  
  This is an immediate consequence of the fact that the the
composition \ $(\TPi)\circ \tj$ coincides with $i$ by
\lcite{\PiTimesRep}, and that $i$  is an injective map.
  \proofend

In particular, this leads to two examples of covariant representations
which are faithfull on $A$, namely the Toeplitz covariant
representation $(\tj,\ts)$ of \lcite{\ToeplitzCovRep} and the regular
covariant representation 
$(i,\n)$ of \lcite{\BundleCovarRep}.

Our next immediate goal will be to show that $(i,\n)$ is strongly
covariant.  In preparation for this we need a technical result to be
proved below.

\state Lemma  
  \label ZinHX
  If $\wdaux$ is any word in $G$ and $X\subseteq G$
is any finite set such that $\spr(\wdaux\inv)\subseteq X$, then
  $$
  \TPi(\tZ_\wdaux)\subseteq \LL{\prd\wdaux}X. 
  $$
  
\proof
Suppose that $\wdaux=(h_1,\cdots,h_n)$ and
let $\gamma=(h_2,\ldots,h_n)$, so that $\wdaux=h_1\gamma$.  Clearly
$\tZ_\wdaux = \tZ_{h_1}\tZ_\gamma$ (closed linear span of products).
Notice that
  $
  \spr(\gamma\inv)\subseteq \spr(\wdaux\inv)\subseteq X,
  $
  so by induction we have that
  $
  \TPi(\tZ_\gamma)\subseteq \LL{\prd\gamma}X.
  $
  This means that the elements of $\TPi(\tZ_\gamma)$ are operators
between the following spaces:
  $$
  \H_X \longrightarrow \H_{\prd\gamma X}.
  $$

  By hypothesis we have that
  $$
  \prd{\gamma\inv}h_1\inv = h_n\inv h_{n-1}\inv \cdots h_2\inv h_1\inv \in
  \spr(\wdaux\inv)\subseteq X.
  $$
  Therefore $h_1\inv\in \prd\gamma X$.  It is also clear that
$\prd{\gamma\inv}\in X$ so that $1\in \prd{\gamma} X$.  In other words
$\{1,h_1\inv\}\subseteq \prd\gamma X$, so $\n_{h_1}\in
\LL{h_1}{\prd\gamma X}$, which means that $\n_{h_1}$ defines an
operator
  $$
  \n_{h_1}: \H_{\prd\gamma X} \longrightarrow \H_{h_1\prd\gamma X} =
\H_{\prd\wdaux X}.
  $$
  The assertion  then follows from \lcite{\Composition}.
  \proofend

\state Proposition
  \label RegIsCovar
  The regular covariant representation  $(i,\n)$ is strongly covariant.

  \proof
  Let $\wd=(g_1,\ldots,g_n)$ be a word in $G$ and let $k\in\tK_\wd$ be
an $\wd$-redundancy.  Recall that $\tK_\wd$ is the closure of the sum
of the $\tZ_\wdaux$, for all words $\wdaux$ such that $\prd\wdaux=1$
and $\spr(\wdaux\inv)\subseteq \spr(\wd)$.  

Since each $\ts_g$ is mapped under $\TPi$ to $B_g$ it is easy to see
that the fact that $\prd\wdaux=1$ implies that $\TPi$ maps $\tZ_\wdaux$
to $B_1$.
  Recall moreover that $B_1$ is defined as the direct limit of the  $\LL
1X$, as $X$ ranges in the collection of all finite subsets of $G$
containing $1$.

We claim that there exists a single $X$ such $\TPi(\tZ_\wdaux)\subseteq
\LL1X$, for all $\wdaux$'s above.  In fact, if $X$ is any finite
subset of $G$ with $\spr(\wd)\subseteq X$ then for every $\wdaux$ with
$\prd\wdaux=1$ and $\spr(\wdaux\inv)\subseteq \spr(\wd)\subseteq X$,
we have by \lcite{\ZinHX} that
  $$
  \TPi(\tZ_\wdaux)\subseteq \LL{\prd\wdaux}X = \LL1X.
  $$

By definition of $\tK_\wd$ we
conclude that
  $
  \TPi(\tK_\wd)\subseteq \LL1X.
  $
  We may therefore regard $\TPi(k)$ as an operator on $\H_X$.

On the other hand we would like to show that $\n_\wd \in \LL
{\prd\wd}{\prd{\wd}\inv X}$.  Since $\ts_\wd\in\tZ_\wd$ and
$\n_\wd=\TPi(\ts_\wd)$ it is enough to prove that 
  $$
  \TPi(\tZ_\wd)\subseteq \LL {\prd\wd}{\prd{\wd}\inv X}.
  $$
  But this  follows once more from Lemma \lcite{\ZinHX} since
  $$
  \spr(\wd\inv) \ \={(\SprOfInverse.i)}\ 
  \prd\wd\inv\spr(\wd) \subseteq
  \prd{\wd}\inv X.
  $$
  Given that $k$ is a redundancy, we have that
  $$
  kb\ts_\wd=0
  \for b\in A,
  $$
  which tells us that $\TPi(k)b\n_\wd=0$.  By what was said above
$\TPi(k)b\n_\wd$ may be interpreted as the composition of operators
  $$
  \H_{\prd{\wd}\inv X} \labelarrow{b\n_\wd}
  \H_X \labelarrow{\TPi(k)} \H_X.
  $$
  Regarding the unit of $A$ as an element of    $\H_{\prd{\wd}\inv
X}$, we therefore have that
  $$
  0 = \TPi(k)b\n_\wd( 1) = 
  \TPi(k)(b).
  $$
  Since this holds for every $b$ in $A$, and since $A$ is dense in
$\H_X$, we conclude that $\TPi(k)=0$.
  \proofend

\state Corollary
  The canonical map 
  $$
  \cpj :A\to \cp
  $$
  (see \lcite{\DefineCp})
  is an injection.

  \proof 
  By \lcite{\DefineCp}  we have that $\cpj = q\circ \tj$.
  Since $(i,\n)$ is strongly covariant by \lcite{\RegIsCovar} we
obtain the *-homomorphism
  $$
  \CpPi : \cp \to C^*(\B)
  $$
  which satisfies
  $
  (\CpPi) \circ j = i,
  $
  by \lcite{\IntroduceStrongly}.
  Since $i$ is injective it follows that $\cpj$  is injective as well.
  \proofend 

We should notice that this does not solve question (1) mentioned after
\scite{\interaction}{7.12}.  The reason is that a single interaction
is not necessarily part of an interaction group.

  \section {Faithful representations}
  As seen in the paragraph following \lcite{\AEmbedsInToep}, strongly
covariant representations $(\pi,\letrep)$, where $\pi$ is faithful,
always exist. However there is no reason for $\rxv\pi\letrep$ to
be faithful.  In this section we will look into the question
of faithfulness of $\rxv\pi\letrep$ carefully.

We therefore fix, throughout this section, a covariant representation
$(\pi,\letrep)$ of the interaction group $(A,G,\m)$ in a C*-algebra
$B$, such that $\pi$ is faithfull.  We will not yet suppose that
$(\pi,\letrep)$ is strongly covariant.

Given that $\pi$ establishes a one-to-one correspondence between $A$
and its image in $B$ there is no loss in generality in assuming that
$A$ is in fact a subalgebra of $B$ and that $\pi$ is the inclusion of
$A$ in $B$.

\definition
  ({\it cf.~\scite{\interaction}{3.6}})
  The covariant representation $(\pi,\letrep)$ will be called
non-degenerated if for every word $\wd$ in $G$ the map
  $$
  a\in A \mapsto a\rep\wd\in B
  $$
  is injective.

For example we have:

\state Proposition 
  \label TwoNonDegenerated
  Both the regular covariant representation $(i,\n)$ and the Toeplitz
covariant representation $(\tj,\ts)$ are non-degenerated.
 
  \proof
  We begin by considering $(i,\n)$.  Let us therefore be given a word
$\wd$ in $G$ and $a$ in $A$ such that $a\n_\wd=0$.  Let $X$ be any
finite subset of $X$ such that $\spr(\wd\inv)\subseteq X$.  By Lemma
$\lcite{\ZinHX}$ we have that $\n_\wd\in\LL{\prd\wd}X.$

Thus, if the unit 1 of $A$ is interpreted as an element of $\H_X$ we
have that
  $$
  0 = a\n_\wd(1) = a.
  $$
  Given that $A$ embedds faithfuly in each $\H_X$ we conclude that
$a=0$.

  In order to show that $(\tj,\ts)$ is  non-degenerated as well it
suffices to observe that, if $a\ts_\wd=0$, then 
  $$
  0 = \CpPi(a\ts_\wd) = a \n_\wd,
  $$
  which implies that $a=0$, by the first part.
  \proofend

Let us write, as usual, $\idp\wd = \rep\wd\rep{\wd\inv}$, where $\wd$
is any  word in $G$.  Since $\idp\wd \rep\wd= \rep\wd$ by
\lcite{\BasicPrepForWords.i}, we have that $a\rep\wd=0$, if and only if
$a\idp\wd=0$.  So the definition of non-degenerated covariant
representations is unchanged if instead we assumed that the map 
  $$
  a \mapsto a\idp\wd
  \eqno{(\seqnumbering)}
  \label InclusionOfRwd
  $$
  is injective.

Recall from  \lcite{\HomosForIdempotnts.ii} that the restriction of
the above map to $\R_\wd$ is a *-homomorphism.  Since injective
*-homomorphisms are necessarily isometric we have:

\state Proposition If $(\pi,\letrep)$ is a non-degenerated covariant
representation and $\pi$ is faithfull then for every word $\wd$ in
$G$ and for every $a\in \R_\wd$ we have that $\|a\idp\wd\| = \|a\| =
\|a\rep\wd\|$.

\proof The argument above clearly gives the first equality.  Moreover
we have
  $$
  \|a\rep\wd\| \leq \|a\| \|\rep\wd\| \leq
  \|a \| =
  \|a\idp\wd \| = 
  \|a\rep\wd\rep{\wd\inv} \| \leq
  \|a\rep\wd\|\|\rep{\wd\inv} \| \leq
  \|a\rep\wd\|.
  \proofend
  $$

A very important feature of non-degenerated covariant representations
is described next:

\state Proposition
  \label NormMwd
  Let $(\pi,\letrep)$ be a non-degenerated covariant representation
of the interaction group  $(A,G,\m)$ in a C*-algebra $B$ such that $\pi$ is faithful.  Then
$\trxv\pi\letrep$ is isometric on $\tM_\wd$ for every word $\wd$ in $G$.

  \proof
  Let $m$ be an element of $\tM_\wd$ of the form
  $$
  m =\sum_{i=1}^n a_i^*\ts_\wd b_i,
  $$
  where $a_1,\ldots,a_n, b_1,\ldots,b_n$ are in $A$. 
We have
  $$
  \|m\|^2 = 
  \Big\| \sum_{i,j=1}^n b_i^*\ts_{\wd\inv} a_i a_j^*\ts_\wd b_j \Big\| =
  \Big\| \sum_{i,j=1}^n b_i^*\m_{\wd\inv}(a_i a_j^*)\tidp{\wd\inv}
b_j \Big\|,
  $$
  where $\tidp\wd = \ts_\wd \ts_{\wd\inv}$, as usual.
  Since $\m_{\wd\inv}$ is a completely positive map by
\lcite{\CompletelyPos.i} there exists an $n\times n$ matrix
$c=\{c_{i,j}\}_{i,j=1}^n$ over $\R_{\wd\inv}$ such that
  $$
  \m_{\wd\inv}(a_i a_j^*) = \sum_{k=1}^n c_{ki}^*c_{kj},
  $$
  for all $i$ and $j$.  We conclude that
  $$
  \|m\|^2 = 
  \Big\| \sum_{i,j,k=1}^n b_i^*c_{ki}^*\tidp{\wd\inv}c_{kj} b_j \Big\|
=
  \Big\| \sum_{k=1}^n d_k^*\tidp{\wd\inv}d_k \Big\|,
  $$
  where $d_k = \sum_{j=1}^n c_{kj} b_j$.  Considering the column matrix
  $$
  y = \pmatrix{\ts_\wd d_1 \cr \vdots \cr \ts_\wd d_n }
  $$
  observe that the expression displayed above says that
  $\|m\|^2 = \|y^*y\|$, but then $\|m\|^2 = \|yy^*\|$ as well.  
  Observe that the $(i,j)$ entry of $yy^*$ is given by
  $$
  \ts_\wd d_i d_j^* \ts_{\wd\inv} =
  \m_\wd(d_i d_j^*)\ts_\wd \ts_{\wd\inv} =
  \m_\wd(d_i d_j^*)\tidp\wd.
  $$
  If one then considers the matrix algebra version of 
\lcite{\InclusionOfRwd}, namely the necessarily isometric
*-homomorphism
  $$
  \{a_{i,j}\}_{ij} \in M_n(\R_\wd)
  \ \longmapsto \
  \{a_{i,j}\tidp\wd\}_{ij} \in M_n(\toep),
  $$
  we conclude that $\|m\|^2$ coincides with the norm of the $n\times
n$ matrix $\{\m_\wd(d_i d_j^*)\}_{ij}$.  It is important to notice
that the outcome of this computation led us to a value which makes no
reference to covariant representations!

Observing  that 
  $
  \trxv\pi\letrep(m) = \sum_{i=1}^n a_i^*\rep\wd b_i,
  $
  we may adopt the exact same strategy to compute the norm of
$\trxv\pi\letrep(m)$, which will therefore give the same result as above.
This is of course due to the fact that the isometric matrix algebra
homomorphism used above, the crucial aspect everything hinges upon, is
based on the non-degeneracy property shared by the Toeplitz covariant
representation and by $(\pi,\letrep)$.  This concludes the proof.
  \proofend

The above result is essentially the same as of the norm calculation in
\scite{\interaction}{4.2}.

  We thus obtain the following crucial faithfulness result:

\state Proposition
  \label FaithfulOnCOne
  Let $(\pi,\letrep)$ be a non-degenerated strongly covariant
representation of $(A,G,\m)$ in a C*-algebra $B$ such that $\pi$ is
faithful.  Then, considering the grading $\{C_g\}_{g\in G}$ of $\cp$
described in \lcite{\CpGrading}, one has that $\rxv\pi\letrep$ is
isometric on $C_g$, for every $g\in G$.

  \proof
  We shall begin by proving that $\rxv\pi\letrep$ is isometric on
$C_1$.
  Observe that, by definition, 
  $$ 
  C_g = \overline{
  \sum_{\prd\wd=1}
  \Z_\wd}.  
  $$
  Since $\Z_X$ is the closed sum of  the   $\Z_\wd$ for 
all words $\wd$ for which 
  $\prd \wd =1$ and
  $\spr(\wd\inv)\subseteq X$,
  one has that 
  $C_1$ is the inductive limit of the $\Z^X$
  for all finite subsets $X$ of $G$.  
  Recalling that $\Z^X$ is a C*-algebra by \lcite{\IntroduceZX}, it is
therefore enough to prove that $\rxv\pi\letrep$ is injective, and
hence isometric, on every $\Z^X$. We therefore suppose that $k\in\Z^X$
is such that $\rxv\pi\letrep(k)=0.$

Assuming without loss of generality that $1\in X$, write
$X=\{x_0=1,x_1,\ldots,x_n\}$ and consider the word
$\wd=(g_1,\ldots,g_n)$, where $g_i=x_ix_{i-1}\inv$, for
$i=1,\ldots,n$.  Clearly $\spr(\wd)=X$, so that $\Z^X=\K_\wd$ as in
\lcite{\DefineKwd}.

Regarding the quotient map 
  $q : \toep \to \cp$,
  it is elementary to prove that $q(\tK_\wd)=\K_\wd$.  So, there
exists $\widehat k\in\tK_\wd$ such that $q(\widehat k)=k$.  For every
$\widehat m\in\tM_\wd$ observe that
  $$
  \trxv\pi\letrep(\widehat k \widehat m) =
  \rxv\pi\letrep(q(\widehat k \widehat m)) =
  \rxv\pi\letrep(k q(\widehat m)) =
  \rxv\pi\letrep(k)\   \rxv\pi\letrep(q(\widehat m)) = 0.
  $$
  By \lcite{\KwdMwdInMwd} we have that  $\widehat k \widehat m\in\tM_\wd$ and by 
\lcite{\NormMwd}   $\trxv\pi\letrep$ is isometric on $\tM_\wd$.  It
follows that $\widehat k \widehat m=0$, and hence that $\widehat k$ is
an $\wd$-redundancy.  Therefore $k=q(\widehat k)=0$.
  This proves that $\rxv\pi\letrep$ is isometric on $C_1$.  

Given any
$g\in G$ and $z\in C_g$ notice that $z^*z\in C_1$.  So
  $$
  \|\rxv\pi\letrep(z)\|^2 = 
  \|\rxv\pi\letrep(z)^* \ \rxv\pi\letrep(z)\| =
  \|\rxv\pi\letrep(z^*z)\| =
  \|z^*z\| = 
  \|z\|^2.
  \proofend
  $$

For example we have:

\state Corollary
  \label RegularFaithOnCOne
  Let $(i,\n)$ be the regular covariant representation of
\lcite{\BundleCovarRep}.  Then $\rxv i\n$ is injective on every  $C_g$.

  \proof Follows immediately from \lcite{\FaithfulOnCOne} once we
realize that  $(i,\n)$ is
  non-degenerated \lcite{\TwoNonDegenerated},
  strongly covariant \lcite{\RegIsCovar},
  and   $i$ is faithful.  
  \proofend

We now want to consider the question of faithfulness of
$\rxv\pi\letrep$.  Observe, however that the hypotheses of
\lcite{\FaithfulOnCOne} are not enough to guarantee that
$\rxv\pi\letrep$ is faithfull.

For example, if $\m_g$ is the identity map on $A$ for every $g$ in
$G$, then one could conceive of a covariant representation
$(\pi,\letrep)$ in which $\pi$ is any faithfull representation of $A$
in $\B(H)$ and $\rep g\equiv1$. It is easy to show that this covariant
representation satisfies the hypotheses of \lcite{\FaithfulOnCOne} but
$\rxv\pi\letrep$ is not faithfull.  To see this notice that, on the
one hand, $\rxv\pi\letrep(\cps_g) = \rep g=1$, but $\cps_g\neq1$ for
$g\neq1$.

To see that in fact $\cps_g\neq1$, consider the amplification
$(\pi',\letrep')$ of $(\pi,\letrep)$ described in
\lcite{\Amplification}.  Then
  $\rxv{\pi'}{\letrep'}(\cps_g)=1\*u_g$, so $\cps_g$ could not
possibly be equal to 1.

  We must therefore give up on the hopes that $\rxv\pi\letrep$ be
faithfull so we will consider its amplification
instead.  But first
let us prove a technical result on graded algebras.

  \state Lemma
  \label FellLemma
  Let 
  $A=\overline{\bigoplus_{g\in G}A_g}$ and 
  $B=\overline{\bigoplus_{g\in G}B_g}$ be topologically graded
C*-algebras and let
  $\phi: A \to B$ be a graded *-homomorphism (meaning that
$\phi(A_g)\subseteq B_g$, for all $g$ in $G$) which is injective on
$A_1$.
  If $G$ is amenable then $\phi$ is injective.

  \proof
  Let us view ${\cal A}=\{A_g\}_{g\in G}$  and  ${\cal
B}=\{B_g\}_{g\in G}$  as  Fell bundles
in the obvious way.
  Since $G$ is amenable one has that ${\cal A}$ is an amenable Fell
bundle by \scite{\amena}{4.7}.  Employing \scite{\amena}{4.2} we
conclude that $A$ is isomorphic to the reduced (or full)  cross
sectional C*-algebra $C_r^*({\cal A})$.

From \scite{\amena}{2.9 and 2.12} it follows that there exists a
faithful conditional expectation 
  $$
  E:A\to A_1
  $$
  which vanishes on $A_g$, for all $g\neq 1$.
  Denote by $F:B\to B_1$ the corresponding conditional expectation.
  It is easy to see that for all $a$ in $A$ one has that
  $$
  \phi(E(a)) = F(\phi(a)).
  $$
  If $a\in A$ is such that $\phi(a)=0$, then
  $$
  0 =   F(\phi(a^*a)) = \phi(E(a^*a)).
  $$
  Since $E(a^*a)\in A_1$, and since $\phi$ is injective on
$A_1$ by hypothesis, we conclude that $E(a^*a)=0$, and
hence that $a=0$ because $E$ is faithful. \proofend

  \state Theorem
  \label AmplificationFaithfull
  Assume that $G$ is an amenable group and let $(\pi,\letrep)$ be a
non-degenerated strongly covariant representation of $(A,G,\m)$ in a
C*-algebra $B$ such that $\pi$ is faithful.  Then
$\rxv{\pi'}{\letrep'}$ is injective, where $(\pi',\letrep')$ is the
amplification of $(\pi,\letrep)$.

  \proof
  Both $\cp$ and $B\*C^*(G)$ are graded algebras and it is easy to 
see that $\rxv{\pi'}{\letrep'}$ is a graded *-homomorphism.  

Observe that for every $z\in C_1$ one has that
  $$
  \rxv{\pi'}{\letrep'}(z) =     \rxv\pi\letrep(z)\*1,
  $$
  so $\rxv{\pi'}{\letrep'}$ is injective on $C_1$ by
\lcite{\FaithfulOnCOne}.  The conclusion then follows from \lcite{\FellLemma}.
  \proofend

Regardless of whether or not $G$ is amenable it is possible to prove
that $\cp$ is isomorphic to the full cross sectional C*-algebra of the
Fell bundle ${\cal C}$, while $B\*C^*(G)$ is the full cross sectional
C*-algebra of the trivial Fell bundle $B\times G$.  Moreover it is
easy to prove that $\rxv{\pi'}{\letrep'}$ is isometric on each fiber
of $\cal C$.  This should be enough to prove that
$\rxv{\pi'}{\letrep'}$ is injective but, without the hypothesis that
$G$ is amenable, I have not been able to find a proof to support this
claim!  The missing argument is: does a fiberwise isometric Fell
bundle homomorphism induces an isomorphism between the full cross
sectional C*-algebras?  
  The answer is unfortunately negative (see
e.g.~\scite{\wasserman}{3.2}).

Since $C^*(\B)$ is already graded we also have:

\state Theorem
  \label RegularIsFaithful
  Assume that $G$ is amenable and let $(i,\n)$ be the regular
covariant representation.  Then $\rxv i\n$ is injective.

  \proof Follows from 
  \lcite{\FellLemma} and \lcite{\RegularFaithOnCOne} as in
  \lcite{\AmplificationFaithfull}.
  \proofend

  \section{Invariant states}
  As before we fix a non-degenerated interaction group $(A,G,\m)$.
Assuming the existence of a faithfull invariant state we will give a
very concrete model for $\cp$.

  \definition A state $\phi$ on $A$ is said to be $\m$-invariant if
for every $g$ in $G$ one has that $\phi\circ\m_g = \phi$.

From now on we fix a $\m$-invariant state $\phi$ and we will let
$\pi$ be the GNS representation of $A$ associated with $\phi$.  The
representation space will be denoted by $H$ and the cyclic vector will
be denoted by $\xi$.

  \state Proposition For every $g$ in $G$ there exists a bounded
linear operator  $\rep g$ in $\B(\H)$ such that
  $$
  \rep g\big(\pi(a)\xi\big) = \pi\big(\m_g(a)\big)\xi
  \for a\in A.
  $$
  In addition, for every $g$ in $G$ one has that $\rep g^* =
\rep {g\inv}$.
  
  \proof Observe that for every $a$ in $A$ we have
  $$
  \|\pi\big(\m_g(a)\big)\xi\|^2 =
  \phi\big(\m_g(a^*)\m_g(a)\big) {\buildrel(\Inequal) \over \leq}
  \phi\big(\m_g(a^*a)\big) = 
  \phi(a^*a) =
  \|\pi(a)\xi\|^2.
  $$
  Therefore the correspondence 
  $$
  \pi(a)\xi\mapsto \pi\big(\m_g(a)\big)\xi  
  $$
  extends to a bounded linear map $\rep g$ on $H$.

  In order to prove the last assertion
  observe that, for every $a,b\in A$, one has
  $$
  \<\rep g\big(\pi(a)\xi\big), \pi(b)\xi\> =
  \<\pi\big(\m_g(a)\big)\xi, \pi(b)\xi\> =
  \phi\big(b^*\m_g(a)\big) \$=
  \phi\Big(\m_g\m_{g\inv}\big(b^*\m_g(a)\big)\Big) =
  \phi\Big(\m_g\m_{g\inv}(b^*)\ \m_g\m_{g\inv}\m_g(a)\Big) \$=
  \phi\Big(\m_g\m_{g\inv}(b^*)\ \m_g(a)\Big) =
  \phi\Big(\m_g\big(\m_{g\inv}(b^*)a\big)\Big) =
  \phi\big(\m_{g\inv}(b^*)a\big) \$=
  \<\pi(a)\xi, \pi\big(\m_{g\inv}(b)\big)\xi\> =
  \<\pi(a)\xi, \rep {g\inv}\big(\pi(b)\xi\big)\>,
  $$
  therefore proving that $\rep g^* = \rep {g\inv}$.
  \proofend

We should remark that $\rep g$ fixes the cyclic vector $\xi$ for all $g$
since
  $$
  \rep g(\xi) =
  \rep g\big(\pi(1)\xi\big) =
  \pi\big(\m_g(1)\big)\xi =
  \pi(1)\xi =
  \xi.
  \eqno{(\seqnumbering)}
  \label ugFixesXi
  $$

  \state Proposition
  \label GNSStrongCov
  Viewing $\letrep$ as a map from $G$ to $\B(H)$ one has that $\letrep$ is a
*-partial representation.  Moreover the pair $(\pi,\letrep)$ is a strongly
covariant
representation  of $(A,G,\m)$ in $H$.

\proof  For every $g$ and $h$ in $G$ and every $a\in A$ we have
  $$
  \rep {g\inv}\rep g\rep h\big(\pi(a)\xi\big) =
  \pi\big(\m_{g\inv}\m_g\m_h(a)\big)\xi =
  \pi\big(\m_{g\inv}\m_{gh}(a)\big)\xi\$=
  \rep {g\inv}\rep {gh}\big(\pi(a)\xi\big),
  $$
  so that $  \rep {g\inv}\rep g\rep h=  \rep {g\inv}\rep {gh}$.
  This proves that $\letrep$ is in fact a partial representation.  Given
$a,b\in A$ we have
  $$
  \rep g \pi(a) \rep {g\inv} \pi(b) \xi =
  \pi\Big(\m_g\big(a\m_{g\inv}(b)\big)\Big) \xi =
  \pi\Big(\m_g(a)\m_g\big(\m_{g\inv}(b)\big)\Big) \xi \$=
  \pi\big(\m_g(a)\big) \pi\Big(\m_g\big(\m_{g\inv}(b)\big)\Big) \xi
=
  \pi\big(\m_g(a)\big) \rep g\rep {g\inv}\pi(b) \xi.
  $$
  This shows that 
  $
  \rep g \pi(a) \rep {g\inv} =
  \pi\big(\m_g(a)\big) \rep g\rep {g\inv},
  $
  and hence that $(\pi,\letrep)$ is a covariant representation.

  In order to prove that  $(\pi,\letrep)$ is strongly covariant
let $\wd$ be a word in $G$ and let $k$ be an $\wd$-redundancy.
Then, for every $b$ in $A$ we have that $kb\ts_\wd  
=0$, whence
  $$
  0 = \trxv\pi\letrep(kb\ts_\wd)\calcat\xi =
  \trxv\pi\letrep(k)\pi(b)\rep \wd\calcat\xi =
  \trxv\pi\letrep(k)\pi(b)\calcat\xi.
  $$
  Since $\xi$ is cyclic we conclude that $\trxv\pi\letrep(k)=0$.
  \proofend 

We now want to consider the question as to whether $\rxv\pi\letrep$ is an
isomorphism onto its range.  Clearly there is no hope for this to be
true unless $\pi$ is faithful.  We therefore assume, from now on,
that $\phi$ is a
  faithful\fn{One does not need $\phi$ to be faithful in order for
$\pi$ to be faithful.  It would be enough to assume that $\phi$
satisfies $(\forall x,y\in A: \phi(xay) = 0 )\Rightarrow a =0.$
However we will soon need $\phi$ to be faithful for other reasons.  }
state,
  which is to say that, for every $a$ in $A$
  $$
  \phi(a^*a) = 0 \implies a =0.
  $$

Recall that we are always working under the standing hypothesis
according to which $(A,G,\m)$ is a non-degenerated interaction group.
However, even if this was not assumed in advance, the existence of a
faithfull invariant state would imply it:

\state Proposition 
  \label FaithStateImplyND   
  Let $(A,G,\m)$ be an interaction group which we exceptionally do not
suppose to be non-degenerated.  If there exists a faithfull invariant
state $\phi$ then $(A,G,\m)$ must be non-degenerated.

  \proof Let $g\in G$ and $a\in A$ be such that $\e_g(a^*a)=0$.  Then
  $$
  0=
  \m_{g\inv}\e_g(a^*a) =
  \m_{g\inv}\m_g\m_{g\inv}(a^*a) =
  \m_{g\inv}(a^*a),
  $$
  so, employing the covariant representation $(\pi,\letrep)$ of
\lcite{\GNSStrongCov} we would have that
  $$
  0 = 
  \<\pi\big(\m_{g\inv}(a^*a)\big)\rep{g\inv}\rep g\xi,\xi\> =
  \<\rep{g\inv} \pi(a^*a)\rep g (\xi),\xi\> \$=
  \<\pi(a^*a)\rep g(\xi),\rep g(\xi)\> =
  \<\pi(a^*a)(\xi),\xi\> =
  \phi(a^*a),
  $$
  and the faithfullness of $\phi$ implies that $a=0$.
  \proofend

Another non-degeneracy property which follows from the existence of a
faithfull invariant state is as follows:

\state Proposition
  \label GNSNonDeg
  If $\phi$ is faithful then $(\pi,\letrep)$ is non-degenerated.

  \proof
  Let $\wd$ be a word in $G$ and let $a\in A$ be such that
$\pi(a)\rep \wd=0$.
Then, 
  $$
  0 = 
  \|\pi(a)\rep \wd \xi\|^2 = 
  \|\pi(a)\xi\|^2 = 
  \<\pi(a)\xi,\pi(a)\xi\> =
  \phi(a^*a),
  $$
  which implies that $a=0$.
  \proofend

The following important result shows that the abstractly defined
crossed product algebra has a very concrete description as an algebra
of operators in the presence of a faithful invariant state:

\state Theorem
  \label ConcreteCp
  Let $(A,G,\m)$ be an interaction group and let $\phi$ be a faithful
$\m$-invariant state.  Let $(\pi,\letrep)$ be the strongly covariant
representation obtained from $\phi$ as in \lcite{\GNSStrongCov} and
let $(\pi',\letrep')$ be its amplification.  If $G$ is
amenable then $\rxv{\pi'}{\letrep'}$ is injective.
  Therefore $\cp$ is isomorphic to the closed *-subalgebra of
$\B(H\*\ell_2(G))$ generated by
  $$
  \{\pi(a)\*1: a\in A\} \cup \{\rep g\*\lambda_g: g\in G\},
  $$
  where $\lambda$ is the left regular representation of $G$.

  \proof
  That $\rxv{\pi'}{\letrep'}$ is injective follows from
\lcite{\AmplificationFaithfull} since $(\pi,\letrep)$ is strongly
covariant \lcite{\GNSStrongCov} and non-degenerated
\lcite{\GNSNonDeg}, and since $\pi$ is obviously faithful.

Moreover, since $G$ is amenable, and hence $C^*(G)$ is nuclear,
$\B(H)\*C^*(G)$ embedds faithfully in $\B(H\*\ell_2(G))$ in such a way
that
  $T\*u_g$ is mapped to $T\*\lambda_g$, for all $T\in\B(H)$ and $g\in G$.
  \proofend

\section{Semigroups of endomorphisms}
  In this section we will discuss the relationship between our notion
of interaction groups and semigroups of endomorphisms.  

We will suppose throughout that $\Pos$ is a subsemigroup of the group
$G$.  In order to simplify certain technical points we will suppose
that $G=\Pos\inv\Pos$.

We will also suppose that $A$ is a unital C*-algebra and that $\ac$ is
a semigroup homomorphism from $P$ to the semigroup of all unital
*-endomorphisms of $A$.  Some will also call $\ac$ an action by
endomorphisms of $\Pos$ on $A$.

Several proposals have appeared in the literature for the notion of
crossed-product of $A$ by $P$ under $\ac$.

  Under the assumption that there exists a faithfull $\ac$-invariant
state $\phi$ on $A$, an assumption that will be enforced henceforth,
the following is perhaps another reasonable proposal: let
$(\pi,H,\xi)$ be the GNS representation of $A$ relative to $\phi$
and, for each $g\in\Pos$, let $\rep g$ be the unique isometry on $H$
such that
  $$
  \rep g\big(\pi(a)\xi\big) =  \pi\big(\ac_g(a)\big)\xi
  \for a\in A.
  \eqno{(\seqnumbering)} 
  \label IsoForSgp
  $$
  One could then define the crossed-product of $A$ by $P$ under $\ac$ to
be the closed *-subalgebra $\sgcp$ of $\B(H)\*C^*(G)$ generated by 
  $$
  \{\pi(a)\*1: a\in A\} \cup \{\rep g\*u_g: g\in \Pos\},
  $$
  where 
  $
  u: G \to C^*(G)
  $
  is the universal representation of $G$.

It is not yet clear to what extent does $\sgcp$ depend on the
faithful invariant state $\phi$ but, based on \scite{\vershik}{6.1}
this is perhaps a sensible definition.

Let us now make the main assumption relating our theory with the
theory of endomorphism crossed-products: 

\sysstate{Hypothesis}{\rm}{\label EndosVsInteraction
  We will suppose, from now on, that there exists an interaction group
$(A,G,\m)$, which also leaves $\phi$ invariant, and such that
  $$
  \m_g=\ac_g
  \for g\in\Pos.
  $$
}

Observe that $(A,G,\m)$ must necessarily be non-degenerated by
\lcite{\FaithStateImplyND}.

The following result says, among other things that, if\/ $\m$ exists,
it can somehow be dug out from the triple $(A,\ac,\phi)$:

  \state Theorem
  \label UniqueExtAndCp
  If\/ $\m'$ is another interaction group satisfying
\lcite{\EndosVsInteraction}, then $\m=\m'$.  In addition, if $G$ is
amenable, then $\cp$ is isomorphic to $\sgcp$, as defined above.

  \proof 
  Let $(\pi,\letrep)$ and $(\pi,\letrep')$ be the covariant
representations on $H$ relative to $\m$ and $\m'$, respectively, as in
\lcite{\GNSStrongCov}.

  Obviously $\rep g = \rep g'$, for all $g$ in $P$, as these agree
with the $\rep g$ defined in \lcite{\IsoForSgp}. 
Since
$\rep{g\inv}=\rep g^*$, we also conclude that $\rep g = \rep g'$, for
all $g$ in $P\inv$.  

For every $g$ in $\Pos$ we know that $\rep g$ is an isometry, and
hence left-invertible.  Se we conclude by \lcite{\RepsCoincidingOnPos}
that $\rep g=\rep g'$, for all $g$ in $G$.
  Next notice that, for all $g\in G$, and $a\in A$, 
  $$
\matrix{
  \rep g \pi(a) \rep {g\inv}\rep g & = &
  \pi\big(\m_g(a)\big)\rep g\rep {g\inv} \rep g =
  \pi\big(\m_g(a)\big)\rep g\cr\cr
  & = & \pi\big(\m'_g(a)\big)\rep g\rep {g\inv} \rep g =
  \pi\big(\m'_g(a)\big)\rep g.
  }
  $$
  It follows that   $\pi\big(\m'_g(a)-\m_g(a)\big)\rep g =0$.  But
since $(\pi,\letrep)$ is non-degenerated by \lcite{\GNSNonDeg}, we
have that $\m'_g(a)=\m_g(a)$.

Refering to the last part of the statement, and using
\lcite{\ConcreteCp}, we must compare the subalgebras of
$\B(H)\*C^*(G)$ generated, on the one hand by
  $$
  \{\pi(a)\*1: a\in A\} \cup \{\rep g\*\lambda_g: g\in G\},
  $$
  and, on the other by
  $$
  \{\pi(a)\*1: a\in A\} \cup \{\rep g\*\lambda_g: g\in \Pos\}.
  $$
  Observe that we are identifying $C^*(G)$ with the reduced C*-algebra
of $G$, and the universal representation $u$ with the regular
representation $\lambda$, since $G$ is supposed to be amenable.

Writing $w_g:= \rep g\*\lambda_g$, 
the problem boils down to whether or not 
all of the  $w_g$ may be generated by those with $g\in \Pos$.

Observing that $w$ is clearly a partial representation, and that $w_g$
is an isometry for all $g$ in $\Pos$, the conclusion follows
immediately from the first part of  \lcite{\RepsCoincidingOnPos},
since we are assuming that $G=\Pos\inv\Pos$.
  \proofend

As already observed, the result above indicates that, in the presence
of a faithful invariant state, there is at most one way to extend the
semigroup action to an interaction group leaving the state invariant.

The extension question is also relevant in the absence of invariant
states:

\sysstate{Question}{\rm}{\label ExtensionQuestion
  Given a semigroup action $\ac$ of $\Pos$ on a C*-algebra $A$, when
is there an interaction group $\m$ such that $\m_g=\ac_g$, for all $g$
in $\Pos$?}

To put matters in perspective let us think of the simplest possible
group-subsemigroup pair, that is, $(\Integers,\N)$.  To give an action
of\/ $\N$ on $A$ is clearly equivalent to giving a single endomorphism
$\ac$ of $A$, which we will suppose injective and unit preserving for
simplicity.  In this case the action may be clearly rerecovered by
iterating $\ac$.

  \state Proposition 
  Let $\ac$ be a unital injective endomorphism of $A$.  
  \izitem
  \zitem If there exists an interaction group $\m$ such that
$\m_1=\ac$, then $\m_1\m_{-1}$ is a conditional expectation onto the
range of $\ac$.
  \zitem If $E$ is a conditional expectation onto the range of $\ac$,
then the operator
  $\L:=\ac\inv\compos E$ is a transfer operator for $\ac$
\scite{\endo}{2.1}, and 
  $$
  \m_n = \left\{\matrix{ \ac^n, & \hbox{if } n\geq0, \cr\cr
                   \L^{-n}, & \hbox{if } n<0.}\right.
  $$
  defines an interaction group such that $\m_1=\ac$.
  \zitem Under the hypothesis of (i) one has that $\m$ is necessarily
given by the expression in (ii) with $E=\m_1\m_{-1}$.

  \proof (i) follows directly from \lcite{\CompletelyPos}.
  Given $E$ as in (ii) one has by \scite{\endo}{2.6} that
$\L:=\ac\inv\compos E$ is a conditional expectation onto the range of
$\ac$.  Observe that it makes sense of speaking of $\ac\inv$ here
since $\ac$ is a bijection onto its range, and since the range of $E$
coincides with that of $\ac$.

  Defining $\m$ as in the statement one can now prove (a tedious but
easy task) that $\m$ is an interaction group.  Clearly $\m_1=\ac$.

Referring to (iii) let $\m$ be given by (ii) and let $\m'$ be another
interaction group such that $\m'_1=\ac$, and
  $$
  \m'_1\m'_{-1} = E = \m_1\m_{-1}.
  $$
  We then have that 
  $$
  \ac\m'_{-1}\m'_1  =
  \m'_1\m'_{-1}\m'_1 =
  \m'_1 = \ac.
  $$
  Given that $\ac$ is injective we conclude that $\m'_{-1}\m'_1$ is
the identity.
  With this observe that
  $$
  \m'_{-1} = 
  \m'_{-1} \m'_1  \m'_{-1} = 
  \m'_{-1} E =
  \m'_{-1} \m_1  \m_{-1} =   
  \m'_{-1} \m'_1  \m_{-1} =   
  \m_{-1}.
  $$
  We therefore conclude that $\m'_g=\m_g$ for all $g\in\{1,-1\}$.
  Since $\m'_1$ is left-invertible and $\m'_{-1}$ is right-invertible,
we can use
\lcite{\OnOneSidedInverse.iv} to 
prove by induction that 
  $$
  \m'_n = (\m'_1)^n
  \for n\geq0,
  $$
  and similarly, using \lcite{\OnOneSidedInverse.iii}, deduce that
  $$
  \m'_n = (\m'_{-1})^{-n}
  \for n<0.
  $$
  Therefore $\m=\m'$.
  \proofend

Thus we see that the extension problem for $(\Integers,\N)$ has a
solution if and only if there exists a conditional expectation onto
the range of $\ac$, and that the collection of all possible solutions
is parametrized by these conditional expectations.

\section{Larsen's crossed products}
  If we consider more general group-subsemigroup pairs $(G,\Pos)$ the
situation may become a lot more complicated.  Let us consider, for
example, the recent work by Larsen \cite{\larsen} on crossed products by
abelian semigroups via transfer operators in which the initial data is
a triple $(\Pos,\ac,\l)$, where $\Pos$ is an abelian semigroup, $\ac$
is an action of $\Pos$ on the (non necessarily unital) C*-algebra $A$,
and $\l$ is an \stress{action by transfer operators}.

We would actually like to consider a slightly different situation,
more general in some aspects, and less general in others: we will
suppose that $A$ is a unital C*-algebra (with this we wish to avoid
the question of extendibility treated by Larsen), $P$ is a
subsemigroup of the non necessarily abelian group $G$ such that
$1\in\Pos$, and $\ac$ is an action of $\Pos$ on $A$ by means of unital
injective endomorphisms.  Observe that $\ac_1$ is necessarily the
identity endomorphism.

We will moreover suppose that we are
given a map (action by transfer operators)
  $$
  \l: \Pos \to \B(A)
  $$
  such that for every $g$ in $\Pos$, $\l_g$ is a transfer operator for
$\ac_g$, and
  $$
  \l_g \l_h = \l_{hg}
  \for g,h\in \Pos.
  $$

Recall that to say that $\l_g$ is a transfer operator for $\ac_g$
means that $\l_g$ is a positive operator on $A$ such that
  $$
  \l_g(a\ac_g(b)) = \l_g(a) b
  \for a,b\in G.
  \eqno{(\seqnumbering)} 
  \label TransferEqn
  $$
  In order to avoid trivialities (such as $\l_g\equiv0$) we will
assume that 
  $$
  \l_g(1)=1
  \for g\in G.
  $$
  Plugging $a=1$ in \lcite{\TransferEqn} then implies that
  $$
  \l_g\circ\ac_g = id_A.
  $$
  In particular $\l_1=id_A$, as well.
  
Question \lcite{\ExtensionQuestion} may then be modified to account
for the $\l_g$ as follows:

\sysstate{Question}{\rm}{\label SecondExtensionQuestion
  Given $(\Pos,\ac,\l)$ as above, is there an interaction group $\m$
such that 
  $\ac_g=\m_g$, and 
  $\l_g=\m_{g\inv}$, for all $g$ in $\Pos$?}

The following is a partial answer:

\state Proposition 
  \label SecondExtensionAnswer
  Suppose that $G=\Pos\inv\Pos$.  Then the above question has an
afirmative answer if and only if
  $
  \ac_g \l_g
  $
  commutes with
  $
  \ac_h \l_h,
  $ 
  for every $g,h\in\Pos$.

  \proof In case $\m$ is an interaction group satisfying the
requirements of \lcite{\SecondExtensionQuestion} we have for all $g$
in $\Pos$ that
  $$
  \e_g:=  \m_g\m_{g\inv} = \ac_g \l_g,
  $$
  so the conclusion follows from \lcite{\BasicPreps.iii}.  Conversely,
given $g$ in $G$ write $g=x\inv y$, with $x,y\in\Pos$, and define
  $$
  \m_g = \l_x\ac_y.
  $$
  We claim that $\m_g$ does not depend on the particular choice of $x$
and $y$.  In order to prove this suppose that $g=z\inv w$, with
$z,w\in\Pos$, as well.  We must then prove that $\l_x\ac_y = \l_z\ac_w$.

As a first case let us consider the situation in which $z=ux$, for
some $u\in \Pos$.  Obviously  $w = zg = uxx\inv y = uy$.
So
  $$
  \l_z\ac_w =
  \l_{ux}\ac_{uy}=
  \l_x\l_u\ac_u\ac_y =
  \l_x\ac_y,
  $$
  because $\l_u\ac_u$ is the identity.
  In the general case pick $u,v\in\Pos$ such that $wy\inv=v\inv u$, so
that $uy=vw$.  Observe that $ux =vwy\inv x = vwg\inv =vww\inv z =
vz. $ So, by the special case already proved, we have
  $$
  \l_x\ac_y =
  \l_{ux}\ac_{uy} =
  \l_{vz}\ac_{vw} =
  \l_z\ac_w.
  $$  
  Therefore the above $\m_g$ is well defined.  It is somewhat curious
that this is so even without having assumed that $P\cap P\inv =
\{1\}$!

We next wish to prove that $\m$ is a partial representation of $G$ on
$A$.  Obviously $\m_1= \l_1\ac_1 = id$.

Given $g,h\in G$ write $g=x\inv y$ and $h=z\inv w$, with
$x,y,z,w\in\Pos$.  Pick $u,v\in\Pos$ such that $zy\inv=v\inv u$ and 
notice that one has
  $$
  uy=vz.
  $$
  Replacing
  $(x,y)$ by $(ux,uy)$, and $(z,w)$ by $(vz,vw)$, we may then assume that
$y=z$.
  With an eye on \lcite{\DefPrep.ii} we compute
  $$
  \m_g\m_h\m_{h\inv} = 
  \l_x\underline{\ac_y\l_y}\ \underline{\vrule width 0pt depth 2.8pt
\ac_w\l_w}\ac_y =
  \l_x\underline{\vrule width 0pt depth 2.8pt \ac_w\l_w}\
\underline{\ac_y\l_y}\ac_y =
  \l_x\ac_w\l_w\ac_y =
  \m_{gh}\m_{h\inv}.
  $$
  Speaking of \lcite{\DefPrep.iii} we have
  $$
  \m_{g\inv}\m_g\m_h =
  \l_y\underline{\vrule width 0pt depth 2.8pt \ac_x \l_x}\
\underline{\ac_y\l_y}\ac_w =
  \l_y\underline{\ac_y\l_y}\ \underline{\vrule width 0pt depth 2.8pt
\ac_x \l_x}\ac_w =
  \l_y\ac_x \l_x\ac_w =
  \m_{g\inv}\m_{gh}.
  $$
  
Our hypotheses clearly imply \lcite{\DefGinter.i-ii}.
  In order to prove \lcite{\DefGinter.iii} let $g=x\inv y$, with
$x,y\in\Pos$.  Take $a$ in $A$ and $b$ in the range of $\m_{g\inv}$,
say $b=\m_{g\inv}(c)$, where $c\in A$.
  Then
  $$
  \m_g(ab) =
  \m_g\big(a\m_{g\inv}(c)\big) =
  \l_x\ac_y\big(a\l_y\ac_x(c)\big) =
  \l_x\Big(\ac_y(a)\ \ac_y\l_y\ac_x(c)\Big) \$=
  \l_x\Big(\ac_y(a)\ \ac_y\l_y\ac_x\l_x\ac_x(c)\Big) =
  \l_x\Big(\ac_y(a)\ \ac_x\l_x\ac_y\l_y\ac_x(c)\Big) \$=
  \l_x\big(\ac_y(a)\big) \l_x\ac_y\l_y\ac_x(c) =
  \m_g(a)\m_g\big(\m_{g\inv}(c)\big) =
  \m_g(a)\m_g(b).
  $$
  The case in which $a$, instead of $b$, lies in the range of
$\m_{g\inv}$  follows by taking adjoints.
  \proofend

We refer the reader to the next section for an example which shows
that the hypotheses of \lcite{\SecondExtensionAnswer} are not always
satisfied.

A situation in which the above result applies is for \stress{linearly
ordered groups}:

\state Proposition
  Let $(\Pos,\ac,\l)$ be as above and suppose that
$G=\Pos\inv\cup\Pos$.  Define
  $$
  \m_g = \left\{\matrix{ \ac_g, & \hbox{if } g\in\Pos, \cr\cr
                   \l_{g\inv}, & \hbox{if } g\inv\in\Pos.}\right.
  $$
  Then $\m_g$ is well defined and $\m$ is an interaction group.

  \proof We will deduce everything from
\lcite{\SecondExtensionAnswer}.  Given $g,h\in\Pos$ we shall prove
that
  $
  \ac_g \l_g
  $
  and 
  $
  \ac_h \l_h
  $
  commute.  Suppose, without loss of generality, that $u:=hg\inv\in
\Pos$, so that $h=ug$.  Then
  $$
  \ac_h \l_h = 
  \ac_{gu} \l_{gu} = 
  \ac_g\ac_u\l_u\l_g.
  $$
  Therefore
  $$
  \ac_h \l_h \ac_g \l_g =
  \ac_g\ac_u\l_u\l_g \ac_g \l_g =
  \ac_g\ac_u\l_u\l_g =
  \ac_g \l_g \ac_g\ac_u\l_u\l_g =
  \ac_g \l_g \ac_h\ac_h.
  $$
  This verifies the hypotheses of \lcite{\SecondExtensionAnswer} and
hence the conclusion follows.
  \proofend

  We therefore see that, under the condition that
$G=\Pos\inv\cup\Pos$, the dynamical systems of \cite{\larsen} are
closely related to interaction groups.

The reader acquainted with Larsen's paper is undoubtedly curious as to
what is the precise relationship between the crossed product defined
in \scite{\larsen}{1.3} and our notion of crossed product, whenever
\lcite{\SecondExtensionQuestion} has an affirmative answer.  However,
while there may be some close relationship between $\toep$ and the
Toeplitz algebra defined in \scite{\larsen}{1.1} it seems that our
notion of redundancy is significantly different from that defined in
\scite{\larsen}{1.2}, so it is unlikely that the crossed products will
coincide.

In addition, given the close relationship between the crossed product
construction studied in \cite{\larsen} and the Cuntz-Pimsner algebras
of product systems (see \scite{\larsen}{3.3}), we also feel that there
is a big discrepancy between our notion of strong covariance and the notion
of Cuntz--Pimsner covariant representations given in
\scite{\fowler}{2.5}.  It is therefore unlikely that our crossed
product construction could be derived from a product system.

  
  \def \half{^{1/2}}
  \def \Hzin{{\cal H}}
  \def \Kzin{{\cal K}}
  \def \ezin{e}
  \def \szin{s}
  \def \Ebig{{\cal E}}  
  \def \F{{\cal F}}
  \def \CAR{A}
  \def \mnwed#1#2{#1_1\wedge\ldots\wedge#1_#2}
  \def \mwed#1{\mnwed #1n}
  \def \mewed#1#2{\ezin(#1_1)\wedge \ldots \wedge \ezin (#1_#2)} 

\section{Example}
  In this section we would like to show an example of a dynamical
system $(A,S,\ac,\L)$, as defined by Larsen in \scite{\larsen}{1}, for
which the hypotheses of \lcite{\SecondExtensionAnswer} do not hold.
We thank Vaughan Jones for suggesting that such an example could be
found by the process of second quantization, as described below.

Having failed to find an example in which $A$ is a commutative algebra
I began wondering if it exists!?

  Let $\Hzin$ be a separable infinite dimensional Hilbert space.  The
CAR algebra of $\Hzin$, denoted $\CAR(\Hzin)$, is the universal
C*-algebra generated by the \stress{annihilation operators}
  $\{a(f)\}_{f\in \Hzin}$ subject to the \stress{canonical
anti-commutation relations}:
  \izitem
  \zitem  $a(f)a(g)+a(g)a(f) =0$,
  \zitem  $a(f)a(g)^*+a(g)^*a(f)=\<f,g\>1$,
  \zitem   $a(f+\lambda g)=a(f)+\overline\lambda a(g)$,
  \medskip\noindent for all $f$ and $g$ in $\Hzin$, and
$\lambda\in\C$.  We refer the reader to \cite{\bratteli} for a
detailed treatment.

\bigskip In order to fix our notation let us briefly describe the
standard representation of $\CAR(\Hzin)$ on the Fermi-Fock space.
Denote by $\F(\Hzin)$ the full Fock space on $\Hzin$, namely
  $$
  \F(H) = \C \oplus \Hzin \oplus \Hzin^{\*2} \oplus \ldots\oplus
\Hzin^{\*n}\oplus \ldots
  $$ 
  and let $P_-$ be the projection defined on each  $\Hzin^{\*n}$ by
  $$
  P_-(f_1\*f_2\*\cdots \*f_n) =
  (n!)\inv \sum_\pi \varepsilon_\pi f_{\pi_1}\*f_{\pi_2}\*\cdots \*f_{\pi_n},
  $$
  where the sum ranges over all permutations $\pi$ of the finite set
$\{1,\ldots,n\}$ and $\varepsilon_\pi$ is the sign of $\pi$.
The Fermi-Fock space $\F_-(\Hzin)$ is defined to be the range of
$P_-$.

The representation of $\CAR(\Hzin)$ on $\F_-(\Hzin)$ we referred to is
defined as follows:
  denoting
  $$
  \mwed f =   P_-(f_1\*f_2\*\cdots \*f_n)
  $$
  one has that
  $$
  a^*(f)(\mwed f) =
  (n+1)\half f\wedge \mwed f.
  $$

In order to obtain an expression for  inner-products in
$\F(\Hzin)$ notice that, given\hfill\break  $f_1,\ldots,f_n,g_1,\ldots,g_n\in\Hzin$
we have
  $$
  \<\mwed f,\mwed g\> =
  \<P_-(f_1\*\cdots \*f_n),P_-(g_1\*\cdots \*g_n)\> \$=
  \<P_-(f_1\*\cdots \*f_n),g_1\*\cdots \*g_n\> =
  (n!)\inv \sum_\pi \varepsilon_\pi 
\<f_{\pi_1}\*\cdots \*f_{\pi_n},g_1\*\cdots \*g_n\> \$=
    (n!)\inv \sum_\pi \varepsilon_\pi 
  \<f_{\pi_1},g_1\>  \ldots \<f_{\pi_n},g_n\>.
  $$
  Thus, if we denote by $\<f,g\>$ the $n\times n$ complex matrix given
by $\<f,g\>_{ij} = \<f_i,g_j\>$, we see that
  $$
  \<\mwed f,\mwed g\> =
  (n!)\inv \det(\<f,g\>).
  $$

  From now on we fix a closed subspace $\Kzin$ of $\Hzin$.  We will
view $\F_-(\Kzin)$ as a closed subspace of $\F_-(\Hzin)$ in the
natural way.

\state Lemma
  Let $\ezin$ denote the orthogonal projection from $\Hzin$ to $\Kzin$ and
let   $\Ebig$ be the orthogonal projection from $\F(\Kzin)$
onto $\F(\Hzin)$.  Then
  $$
  \Ebig(\mwed f) = \mewed f n.
  $$

\proof
  Given $f_1,\ldots,f_n\in\Hzin$ we claim that 
  $$
  y:=  \mwed f -
  \mewed f n \in \F(\Kzin)^\perp.
  $$
  To see this let $g_1,\ldots,g_n\in\Kzin$ and notice that
  $$
  \<y,\mwed g\> = 
  \<\mwed f,\mwed g\> -
  \<\mewed f n,\mwed g\> \$=
  { \det(\<f,g\>) - \det(\<\ezin(f),g\>) \over n!} = \cdots
  $$
  Since 
  $
  \<\ezin(f_i),g_j\> = 
  \<f_i,\ezin(g_j)\> = \<f_i,g_j\>,
  $
  we conclude that the above vanishes.
  The conclusion then follows easily.
  \proofend

\state Lemma
  \label EquationWithE
  If $f\in\Hzin$ then 
  $
  a(f)\Ebig = \Ebig a(\ezin(f)) = a(\ezin(f)) \Ebig.
  $

  \proof We prove instead that
  $\Ebig a^*(f)  = a^*(\ezin(f)) \Ebig= \Ebig a^*(\ezin(f)).$
  Given $f_1,\ldots,f_n\in\Hzin$ we have
  $$
  \Ebig a^*(f)(\mwed f) =
  (n+1)\half   \Ebig (f\wedge\mwed f) = 
  (n+1)\half   \ezin(f)\wedge \mewed fn \$= 
  a^*(\ezin(f)) \mewed fn = 
  a^*(\ezin(f)) \Ebig(\mwed f),
  $$
  proving that $\Ebig a^*(f) = a^*(\ezin(f)) \Ebig$.  The remaining
equality follows from this upon replacing $f$ by $\ezin(f)$.
  \proofend

It is interesting to observe that this result implies that
$\CAR(\Kzin)$ commutes with $\Ebig$.

The following is probably a well known result for specialists in
second quantization but we were unable to find a reference for it:

\state Proposition
  \label CarCondExpec
  There exists a unique bounded linear map $E:\CAR(\Hzin)\to \CAR(\Kzin)$
such that
  $$
  E\Big(a^*(f_1)\ldots a^*(f_n)a(g_1)\ldots a(g_m)\Big) = 
  a^*\big(\ezin(f_1)\big)\ldots a^*\big(\ezin(f_n)\big)
  a\big(\ezin(g_1)\big)\ldots a\big(\ezin(g_m)\big),
  $$
  for every 
  $f_1,\ldots, f_n, g_1,\ldots, g_m\in\Hzin$.
Moreover $E$ is a  conditional expectation onto $\CAR(\Kzin)$.

  \proof
  We have by \lcite{\EquationWithE} that 
  $$
  \Ebig a^*(f_1)\ldots a^*(f_n) a(g_1)\ldots a(g_m)\Ebig =
  a^*(\ezin(f_1))\ldots a^*(\ezin(f_n)) \Ebig a(\ezin(g_1))\ldots
a(\ezin(g_m)) \$=
  a^*(\ezin(f_1))\ldots a^*(\ezin(f_n)) a(\ezin(g_1))\ldots a(\ezin(g_m))\Ebig.
  $$
  Since $\CAR(\Kzin)$ is faithfully represented on $\F_-(\Kzin)$, and
since $\F_-(\Kzin)$ is precisely the range of $\Ebig$, we conclude
that there exists a unique bounded linear map $E:\CAR(\Hzin)\to
\CAR(\Kzin)$ such that
  $$
  \Ebig x \Ebig = E(x) \Ebig
  \for x\in \CAR(\Hzin).$$
  It is now easy to verify that $E$
satisfies the required conditions.
  \proofend

Let $\szin$ be an isometry on $\Hzin$.  By the universal property of
$\CAR(\Hzin)$ there exists a unique unital *-endomorphism $\ac$ of
$\CAR(\Hzin)$ such that
  $$
  \ac(a(f)) = a(\szin(f))
  \for f\in \Hzin.
  \eqno{(\seqnumbering)}
  \label CAREndomorphism
  $$

Since $\CAR(\Hzin)$ is a simple algebra $\ac$ is necessarily
injective.  The range of $\ac$ obviously coincides with $\CAR(\Kzin)$,
where $\Kzin$ is the range of $\szin$.

As in \scite{\endo}{2.6} we may produce a transfer operator $\L$ for
$\ac$ by setting
  $$
  \L = \ac\inv \compos E,
  \eqno{(\seqnumbering)}
  \label CARTransfer
  $$
  where $E$ is as in \lcite{\CarCondExpec}.  It is elementary to
verify that 
  $$
  \L\Big(a^*(f_1)\ldots a^*(f_n)a(g_1)\ldots a(g_m)\Big) = 
  a^*\big(\szin^*(f_1)\big)\ldots a^*\big(\szin^*(f_n)\big)
  a\big(\szin^*(g_1)\big)\ldots a\big(\szin^*(g_m)\big).
  $$

Let us now get a bit more concrete and consider $\Hzin=\ell^2(\N)$.
Let $\szin_1$ and $\szin_2$ be two commuting isometries on $\Hzin$
such that $\szin_1\szin_1^*$ and $\szin_2\szin_2^*$ do not commute.
In order to construct these in a concrete way let $\szin_1$ be the
unilateral shift on $\ell^2(\N)$ and, identifying $\ell^2(\N)$ with
the Hardy space $H^2$ as usual, let $\szin_2=T_\phi$ be the Toeplitz
operator whose symbol is the Blaschke factor
  $$
  \phi(z) = {z - a \over 1-\bar a z}
  \for z\in\C,
  $$
  where $0<|a|<1$.  We leave it for the reader to prove that $\szin_1$
and $\szin_2$ are indeed commuting isometries with noncommuting final
projections.

Let $\ac_1$ and $\ac_2$ be the endomorphisms of $\CAR(\Hzin)$
respectively obtained from $\szin_1$ and $\szin_2$ as in
\lcite{\CAREndomorphism}.  Also let $\L_1$ and $\L_2$ be the
corresponding transfer operators obtained by \lcite{\CARTransfer}.

Observe that the fact that $\szin_1$ and $\szin_2$ commute implies
that $\ac_1$ and $\ac_2$ commute.  
  Obviously $\szin^*_1$ and $\szin^*_2$ commute as well which entails
the commutativity of $\L_1$ and $\L_2$.

For every $(n,m)\in\N\times\N$ set
  $$
  \ac_{(n,m)} = \ac_1^n\ac_2^m
  \and
  \l_{(n,m)} = \L_1^n\L_2^m.
  $$
  It is now elementary to check that $(\CAR(\Hzin),\N\times\N,\ac,\l)$
is a dynamical system as defined by Larsen in \scite{\larsen}{1}.
  We shall prove however that it does not satisfy the hypothesis of
\lcite{\SecondExtensionAnswer}.  In fact, notice that
  $$
  \ac_{(1,0)}\Big(\l_{(1,0)}\big(a(f)\big)\Big) =
a\big(\szin_1\szin_1^*(f)\big),
  $$
  while 
  $$
  \ac_{(0,1)}\Big(\l_{(0,1)}\big(a(f)\big)\Big) =
a\big(\szin_2\szin_2^*(f)\big),
  $$
  from which one immediately sees that 
  $
  \ac_{(1,0)}\l_{(1,0)}
  $
  and 
  $
  \ac_{(0,1)}\l_{(0,1)}
  $
  do not commute.

  Needless to say, there is no interaction group extending $\ac$ and
$\l$, precisely by Theorem \lcite{\SecondExtensionAnswer}.

\section{Appendix}
  In this section we would like to give a detailed proof of the
statement made in the introduction, according to which one cannot
always find an isometry $S$ in $\O_A$ such that
  $\ac(a) = SaS^*$, for all $a$ in $C(K)$, where $K$ is Markov's space
and $\ac$ is the endomorphism of $C(K)$ induced by the Markov
subshift.

To focus on a single counterexample we will let
  $$
  A=\pmatrix{1 & 1\cr 1 & 1}
  $$
  in which case $\O_A$ coincides with $\O_2$.  We will denote by $B$
the standard copy of the CAR algebra inside $\O_2$, and by $D$ the
diagonal subalgebra of $B$, which is isomorphic to $C(K)$.  It is well
known that $D$ is a maximal abelian subalgebra of $B$. 

A crucial tool in our argument below is the fact that $D$ is also a
maximal abelian subalgebra of $\O_2$ (see \cite{\cuntzTwo}).  This
also follows from \scite{\renault}{Proposition II.4.7} given the
groupoid description of $\O_2$ 
\scite{\renault}{Section III.2}.

In order to prove the inexistence of an isometry $S$ in $\O_2$
satisfying \lcite{\FirstCovariance} let us argue by contradiction and
hence suppose that such an $S$ exists.  Since $\ac(1)=1$, one must
have that $S$ is unitary.

Letting $D^+$ be the necessarily abelian subalgebra of $\O_2$ given by
$$D^+=S^*DS,$$ we claim that $D\subseteq D^+$. In fact, given $a\in D$
we have
  $$
  a = S^*SaS^*S = S^*\ac(a)S \in S^*DS = D^+.
  $$
  
Since $D$ is maximal abelian in $\O_2$ we conclude that $D=\D^+$.  For
every $a\in D$ we then have that $b:= S^*aS\in D$ and
  $$
  a = SS^*aSS^* = SbS^* = \ac(b),
  $$
  from which one would deduce the absurdity that $\ac$ is surjective.
Thus, no such $S$ may exist!


  \endgroup

  \begingroup
  \bigskip \bigskip  \font \sc = cmcsc8 \sc
  \parindent = 0pt

  Departamento de Matem\'atica

  Universidade Federal de Santa Catarina

  88040-900 -- Florian\'opolis -- Brasil

  \medskip
  {\tt exel@\kern1pt mtm.ufsc.br} 

  \endgroup
  \end